 \documentclass[12pt]{article} 
\usepackage[pdftex]{graphicx} 
\usepackage{amsfonts}
\usepackage{mathrsfs}

\newcommand{\eq}[1]{$(\ref{#1})$}
\newcommand{\R}{\mathbb{R}}
\newcommand{\bS}{{\bf S}}
\newcommand{\CC}{{\cal C}}
\newcommand{\Po}{\mathscr{P}}
\newcommand{\dist}{{\rm dist}}
\renewcommand{\epsilon}{\varepsilon}

\newtheorem{theo}{Theorem}[section]
\newtheorem{prop}{Proposition}[section]
\newtheorem{lemm}{Lemma}[section]
\newcommand{\bea}{\begin{eqnarray}}
\newcommand{\eea}{\end{eqnarray}}
\newcommand{\beaa}{\begin{eqnarray*}}
\newcommand{\eeaa}{\end{eqnarray*}}

\newcommand{\eqco}{\setcounter{equation}{0}}

\title{{Strict inequalities of critical values in continuum percolation}}

\author{Massimo Franceschetti\footnote{Department of Electrical and Computer Engineering,
University of California San Diego,
La Jolla, CA, 92093-0407,
Email: massimo@ece.ucsd.edu}
 \footnote{Partially supported by the Isaac Newton Institute, Cambridge, and by NSF grant CNS0916778.}
Mathew D. Penrose\footnote{ Department of Mathematical Sciences,
 University of Bath BA2 7AY, United Kingdom,
Email: m.d.penrose@bath.ac.uk}
 \footnote{ Partially supported by the Isaac Newton Institute, Cambridge}
 and Tom Rosoman\footnote{ Department of Mathematical Sciences,
 University of Bath BA2 7AY, United Kingdom,
Email: ter20@bath.ac.uk} \footnote{ Supported by an EPSRC studentship}
}

\date{\today}

\begin{document}

\maketitle

\begin{abstract}\addcontentsline{toc}{section}{Abstract}

We consider the supercritical finite-range random connection model 
where the points $x,y$ of a homogeneous planar Poisson process
are connected with probability $f(|y-x|)$ for a given $f$. 
Performing percolation on the resulting graph,
we show that the critical probabilities   for site and bond percolation 
satisfy the strict inequality $p_c^{\rm site} >   p_c^{\rm bond}$. 
We also show that reducing the connection function $f$ strictly   
increases the critical Poisson intensity.

Finally, we deduce that performing  a spreading transformation 
on $f$ (thereby allowing connections over greater distances
but with lower probabilities, leaving average degrees unchanged) 
 {\em strictly} reduces the critical Poisson
intensity.  This is of practical relevance,
 indicating that in many real networks it is in principle possible
 to exploit the presence of spread-out, long range connections, to
 achieve connectivity at a strictly lower density value.
\end{abstract}



Short title: Continuum percolation inequalities

\section{Introduction}
\label{secintro}

Since exact formulae for critical values in percolation are known
only for a few special cases, it is of interest
to obtain partial information in the form
of inequalities between critical values for different 
percolation models. This is especially true for continuum
percolation; no exact critical values at all are known
in the continuum, while on the other hand some interesting
inequalities have been discovered.  A striking result of this
type says that for  percolation of copies of a fixed
convex shape of unit area centred at Poisson points in the plane,
the critical intensity is less for a triangle than for
any other shape. This was established as a weak inequality
by Jonasson (2001), and as  a strict inequality by
Roy and Tanemura (2002). 

The present paper is concerned with
 another result of this type, which says that  for the random connection
model over Poisson  points in the plane (or in higher dimensions),
the critical intensity is decreased under a spreading transformation
of the connection function, whereby connections between
more distant points are allowed but with lower probability, so
that the average degree remains unchanged.

A related topic is the comparison of critical
values for bond and site percolation. Given
an infinite connected graph $G$, let us denote
these critical values by $p_c^{\rm bond}$ and $p_c^{\rm site}$, respectively.
The weak inequality $p_c^{\rm site}\geq  p_c^{\rm bond}$ 
can easily be proven by dynamic coupling, 
see for example Chapter~2 of Franceschetti and
 Meester~(2007).  If $G$ is a rooted tree,
 then it is easy to see that $p_c^{\rm site} =  p_c^{\rm bond}$, as each
vertex, other than the root, can be uniquely identified by an edge and vice versa. By adding finitely many edges to an infinite tree, one can also construct other connected graphs for which the equality holds.  
On the other hand, the strict 
 inequality $p_c^{\rm site} >  p_c^{\rm bond}$ has also been shown to hold in
many circumstances. Grimmett and Stacey~(1998) proved it for
a large class of `finitely transitive' graphs including the
  $d$-dimensional hypercubic lattices. 

In this paper we show
$p_c^{\rm site} >  p_c^{\rm bond}$ 
for certain
 {\em random} graphs arising
in {\em continuum} percolation. Such graphs are not covered
by previous results because they
are not finitely transitive;
  since their node degrees are not bounded,
 the group action defined by their automorphisms has
 infinitely many orbits.
Continuum percolation graphs
are of particular  interest in the context of communication networks and  are treated extensively in the books by  Franceschetti and Meester~(2007), Meester and Roy~(1996), and Penrose~(2003).

We consider the 
 {\em random connection
model} (RCM) of continuum percolation, which is defined as follows.
 Let $\lambda >0$ and let
 $f: \mathbb{R}^{+} \rightarrow [0,1]$ 
(the so-called {\em connection function}) be specified.
Let $\Po_\lambda$ be a homogeneous
 Poisson point process in the plane of intensity $\lambda$ and  connect every pair of points $x,y \in \Po_\lambda$  with probability $f(|x-y|)$,
 where $|\cdot|$ denotes the Euclidean norm 
({\em Gilbert's graph} is the special case of the
 RCM with $f \equiv {\bf 1}_{[0,1]}$).
Provided  $\lambda $ exceeds a critical value $\lambda_f$ which depends
on $f$, the RCM graph has an infinite component almost surely.

\begin{figure}
\begin{center}
\scalebox{.8}{\includegraphics{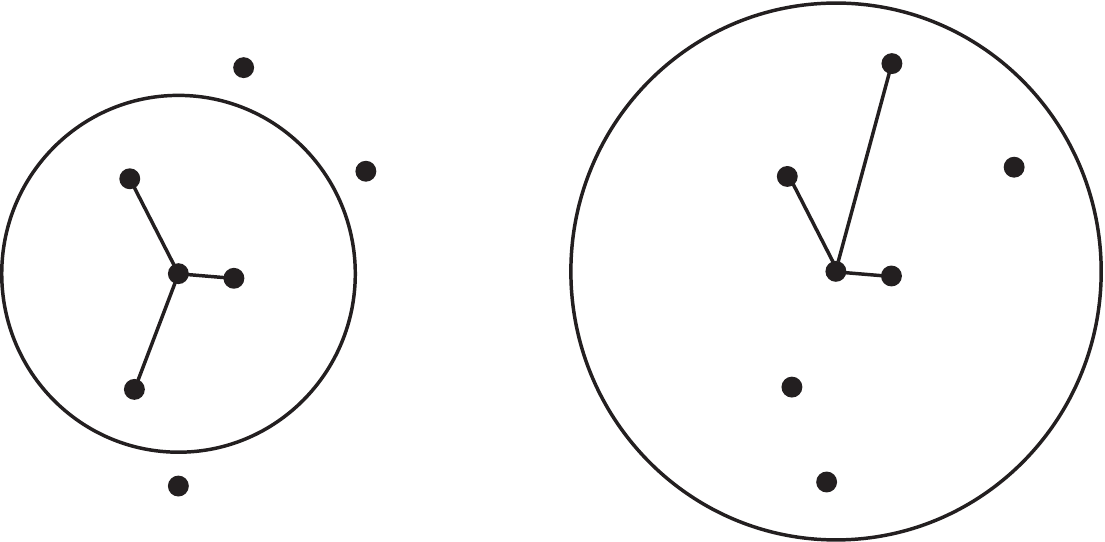}}
\end{center}
\caption{ Left-hand side: RCM with connection function $f(r)$. 
Right-hand side: RCM with connection function $p \cdot f(\sqrt{p}r)$. 
Only connections to the centre point are depicted.
}
\label{spread}
\end{figure}

In Theorems \ref{bool} and \ref{RCM}, we prove
 that $p_c^{\rm site} > p_c^{\rm bond}$ for the random graphs
arising from the supercritical RCM using any nonincreasing
connection function with finite range
 (including Gilbert's graph). In Theorem \ref{th:squash} we
show that
 replacing the connection
function $f(\cdot)$ by a smaller connection function $g(\cdot)$
causes the critical intensity to strictly increase, 
i.e. $\lambda_{g} > \lambda_f$.
  In Theorem \ref{th:spread} we consider
 the spreading transformation already mentioned, which is
defined as follows. Given connection function $f$ and given $0<p<1$,
define the spread-out connection function $S_pf$ by
\begin{equation}
S_pf (r) = p \cdot f(\sqrt{p} r).
\label{Spdef}
\end{equation}
 Thus, the probabilities are reduced by a factor $p$ but the function is
 spatially stretched so as to maintain the same expected number of 
connections per node; see Figure~\ref{spread} for a visual representation.  

Theorem \ref{th:spread} says that
 for any nonincreasing  connection
function $f$ with finite range, and for $0< p <1$,
 we have $\lambda_{S_pf} < \lambda_f$ 
(the weak version of this inequality is much simpler, see
 Franceschetti et al.~(2005)).
Consequently, as $p \downarrow 0$
 the approach of the critical value $\lambda_{S_pf}$ 
to its limiting value (which is known to equal
$(2 \pi \int_0^\infty r f(r)dr )^{-1}$, see  Penrose ~(1993))
is strictly monotone. 


In applications 
 this is of interest, as it shows that unreliable, spread-out connections are strictly advantageous for reaching connectivity at  a given node density value. In communication networks, for example, the `quality' of a communication link decreases as the distance between transmitter and receiver increases. Hence, connections can be established between nearby nodes, but reliable long-range connections are more difficult to obtain.  Our results show that highly reliable, but short-range connections could be exchanged with less reliable, but longer-range connections, to obtain network connectivity more easily (in a strict sense), provided that the average number of functioning connections per node remains the same.

Our results carry over to a large class of connection functions
having {\em infinite} range, i.e. with unbounded support.
However, the proofs for this case require  significant modifications 
and lengthy further arguments,
 and to keep the length of the current paper  
under control, we shall deal with the infinite range case 
elsewhere.

We conclude this section with some observations leading to some
 open problems. First, there exist spreading transformations for
 which the percolation threshold is unaffected. For example, an affine
 transformation of the plane converting discs into aligned ellipses of
 the same area would spread-out connection lengths but would not affect the
 percolation threshold. Second, there are spreading transformations for 
which neither  weak nor strong inequalities are known. For example, the
 effect of the shifting and squeezing transformation for annuli 
 considered by Franceschetti et al.~(2005) and independently by
 Balister, Bollob\'{a}s, and Walters~(2004) is known only in the
 spread-out limit. Similar limiting results as the dimension of the 
space spreads to infinity are given by Meester, Penrose, and Sarkar~(1997).

We are not aware of any lattice analogue of Theorem \ref{th:spread}.
In a lattice version of the spreading transformation, where
one considers e.g. bond percolation on the vertices of
 ${{\mathbb Z}}^d$ with range $r$ becoming large,
the  critical value is known to approach its branching process
limit (see Penrose~(1993) or Bollob\'{a}s, Janson, and Riordan~(2005)),
 but the convergence is not known to be monotone.

\section{Statement and discussion of results}
\eqco

For $\lambda >0$, let  $\Po_\lambda$ be a homogeneous
   Poisson process in $\R^2$ of intensity $\lambda$.
The random connection model (RCM) driven by $\Po_\lambda$,
 with nonincreasing connection function $f: \mathbb{R}^{+} \rightarrow [0,1]$,
is obtained by 
 connecting each 
 pair of points $x,y \in \Po_\lambda$
by an undirected edge
  with probability $f(|x-y|)$, 
independently of other pairs.
We denote the resulting graph by $RCM(\lambda,f)$.
For a formal description of the RCM, see Meester and Roy ~(1996). 
It is well known that provided
$0<\int_0^\infty r f(r) dr<\infty,$
 the RCM has a critical density value $\lambda_f \in (0, \infty)$, in
fact with
\begin{equation}
\lambda_f \geq \frac{1}{2 \pi \int_0^\infty r f(r) dr },
\label{100122a}
\end{equation}
  such that if $\lambda > \lambda_f$ then there exists a.s.\ a unique infinite connected component in $RCM(\lambda,f)$, while  if
 $\lambda < \lambda_f$ then there is a.s no infinite connected component
in $RCM(\lambda,f)$, see Penrose (1991), Meester and Roy (1996).
 When it exists, we denote this infinite component by $\CC$.
 
In the site percolation model on  $\CC$, each vertex is independently marked open with probability $p$, and closed otherwise, and 
we look for an unbounded connected component in  the 
 subgraph  $\CC_v$
induced by the open vertices.  It is easy to see that this is equivalent
 to rescaling the original Poisson process to one with intensity
 $p \lambda$ and looking for an unbounded connected component there. 
It follows that for $\lambda > \lambda_f$ there is a critical
 value $p_c^{\rm site} \in (0,1)$ 
(namely $p_c^{\rm site} = \lambda_f / \lambda$) such that if
 $p > p_c^{\rm site}$ then there is a.s.\ an infinite connected component
 in $\CC_v$, and if $p < p_c^{\rm site}$ then there is
 a.s.\ no such infinite component.

In the  bond percolation model on $\CC$, we independently declare each edge
to be open with probability $p$, and closed otherwise, and look for an
 unbounded connected component in the subgraph $\CC_e$ induced by the open
 edges. This is equivalent to constructing an RCM with the original
 connection function $f$ replaced by $p \cdot f$.
There is a critical probability $p_c^{\rm bond} \in (0,1)$ such that if $p > p_c^{\rm bond}$ then there is a.s.\ an infinite connected component in 
$\CC_e$, and if $p < p_c^{\rm bond}$ then there is 
a.s.\ no such  infinite component. 
To verify  that $p_c^{\rm bond} \in (0,1)$, observe that
$p_c^{\rm bond} \leq p_c^{\rm site} <1$ while by (\ref{100122a}) 
we have $p_c^{\rm bond} \geq (2 \pi \lambda \int_0^\infty r f(r) dr)^{-1} $.


A special case of the RCM arises when $f(r) = 1$ for $r \leq 1$ and $f(r) = 0$ for $r > 1$. This is called
Gilbert's graph, $G(\Po_ \lambda,1)$, and is formed from
$\Po_\lambda$
by joining every two points $x,y \in \Po_\lambda$ with
 $|x - y| \leq 1$. 
 We denote the critical value of $\lambda$ by $\lambda_c$ in this case.

Our first results provide strict inequalities between 
$p_c^{\rm site}$ and $p_c^{\rm bond}$.
We first prove our results for Gilbert's graph
(Theorem \ref{bool}) and then generalise to
the random connection model;
the proof for the first case sets
 up many of the arguments in the more general case.

\begin{theo} \label{bool}
Consider  $G(\Po_ \lambda,1)$ for $\lambda>\lambda_c$.
On $\CC$ we have $p_c^{\rm site} > p_c^{\rm bond}$.
\end{theo}

The next result generalizes Theorem \ref{bool}, and
 concerns the RCM with connection function $f$ having  bounded support.
\begin{theo}\label{RCM}
Consider
  $RCM(\lambda,f)$
 for $\lambda>\lambda_f$. If $f$ is nonincreasing and $0 < \sup\{a:f(a) > 0\} < \infty$
then  on $\CC$ we have
$p_c^{\rm site} > p_c^{\rm bond}$.
\end{theo}

Our next result provides a strict inequality governing
 the effect on the RCM crictial intensity
$\lambda_f$ if one reduces the connection function $(f(r), r \geq 0)$
by a constant factor, i.e. if one uses 
instead the `squashed' connection function $pf := (pf(r), r \geq 0)$ with
$p \in (0,1)$ a constant.
The weak inequality $\lambda_{pf} \geq \lambda_f$ is clear,
and the next result improves it to a strict inequality.
\begin{theo} \label{th:squash}
Let $q_0 \in (0,1)$. Suppose the connection function
 $f$ is nonincreasing and $0 < \sup\{a:f(a) > 0\} < \infty$.
Then $\lambda_{q_0f} > \lambda_f$.
\end{theo}
In fact, Theorem \ref{th:squash} holds in greater generality:
if
$f$ and $g$ are
 connection functions,
 satisfying the hypotheses of Theorem \ref{RCM},
with
 $g(r) \leq f(r)$ for all $r$ and  $g(r) < f(r)$
for $r$ in some sub-interval of $(0,\infty)$, then $\lambda_g
<\lambda f$.  This can be proved by an extension of the 
proof of Theorem \ref{th:squash}, which we omit.

Our last result is concerned with the spreading-out transformation
$S_p$ defined by (\ref{Spdef}).
\begin{theo}
\label{th:spread}
Suppose $ 0 < p < q \leq 1$.
Suppose the connection function
$f$
 is nonincreasing and $0 < \sup\{a:f(a) > 0\} < \infty$.
Then $\lambda_{S_pf} < \lambda_{S_q f}$. Also,
the inequality  \eq{100122a} is strict. 
\end{theo}

To conclude this section we 
 give an overview of the technique of proof and related
literature. We note first that 
our proof of all of these results easily extends to $3$ or more dimensions.

The basic strategy is to adapt the enhancement
 technique developed for percolation on lattices by
 Menshikov~(1987),
 Aizenman and Grimmett~(1991),
 Grimmett and Stacey~(1998).
 This consists of constructing an `enhanced' version of the site
 percolation process  (i.e., one with some extra open sites added according
to certain rules),
for which the critical probability is 
\emph{strictly less} than that of the original site process.
 Then one can use dynamic coupling of the enhanced model with
 bond percolation to complete the proof.

We face two main difficulties when trying to extend the enhancement technique
 to a continuum random setting. One of these  amounts to constructing
 the desired enhancement on a random graph rather than on a deterministic
 one. The second one consists in adapting some basic inequalities for
 the enhanced graph, given in the discrete setting by Aizenman and
 Grimmett~(1991), to the continuum setting.  Because
the possible configurations outside a given region
now provide a continuum of possible boundary conditions, 
this requires somehow more involved geometric constructions
 and a careful incremental 
build-up of the Poisson point process. Once we circumvent these obstacles,
 it is not too difficult to obtain the final result using a classic dynamic 
coupling construction.

 In order to keep the main ideas of the proof clear, we 
first prove Theorem \ref{bool},
 and later adapt the proof
 to the general case of Theorem \ref{RCM}.
 The proof of Theorem \ref{th:squash}
uses an argument involving `diminishment', rather than enhancement,
of site percolation, and the proof of Theorem
\ref{th:spread} uses the preceding results along with 
a coupling argument related to that used by Franceschetti
et al.~(2005) to get the weak version of
Theorem \ref{th:spread}.

The enhancement strategy has proven useful to show strict inequalities in a variety of contexts: Bezuidenhout, Grimmett, and Kesten~(1993), and Grimmett~(1994), use this technique in the context of Potts and random cluster models;  Roy, Sarkar, and White~(1998) use it in the context of directed percolation.
In the continuum, Sarkar~(1997) uses enhancement to demonstrate coexistence
of occupied and vacant phases for the three-dimensional
  Poisson Boolean model.
  Roy and Tanemura~(2002) use it in the context of percolation of different convex shapes.

\section{Gilbert's Graph: Proof of Theorem \ref{bool}}
\label{secgil}
\eqco

\begin{figure}[htbp]
\begin{center}
\scalebox{.8}{\includegraphics[angle = 270]{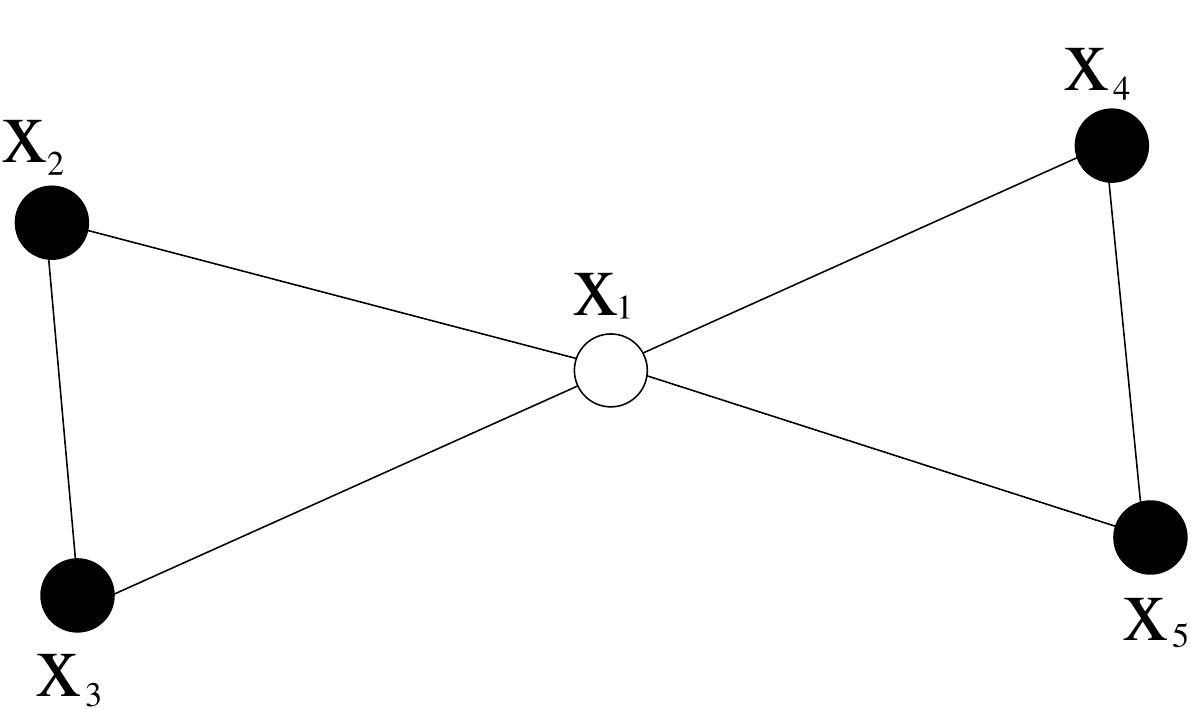}}
\end{center}
\caption{The bow tie enhancement.}
\label{fig:enh1}
\end{figure}

We now describe the enhancement needed to prove Theorem $1$.
 Throughout this section we consider Gilbert's graph 
$G(\Po_ \lambda,1)$ with $\lambda > \lambda_c$.
The
objective is to describe a way to to add open vertices to the site
percolation model 
 without changing the coupled bond percolation model. To do so, we
introduce two kinds of coloured vertices, red vertices (the original
open vertices) and green vertices (closed vertices which have been
enhanced) and for any two vertices $x,y$ we write that $x \sim y$ if
they are joined by an edge. In $G(\Po_\lambda,1)$, if we have vertices
$v,w,x,y,z$ such that $x$ is closed, has no neighbours
other than $v,w,y,z$, which are all red, and $v \sim w$ and
$y \sim z$ but there are no other edges amongst $v, w, y$
and $z$ then we say $x$ is {\em correctly configured}
 in $G(\Po_\lambda,1)$, and refer to this as a {\em bow tie}
 configuration of edges. If a vertex $x$ is correctly configured
we make it green with probability $q$, independently of everything
else; see Figure~\ref{fig:enh1}. 

Let $B_n$ be the open disc of radius $n$ centred at the origin. Let
$\underline{Y} = (Y_i, i \geq 0)$  and $\underline{Z} = (Z_i, i \geq
0)$ be sequences of independent uniform $[0,1]$ random variables.
List the vertices of $\Po_{\lambda}$ in order of increasing distance
from the origin as $x_1,x_2,x_3,\dots$ . Declare a vertex $x_i$ to be 
{\em red} if
$Y_i < p$ and {\em closed} otherwise. Once the sets of red and closed
vertices have been decided in this way, apply the enhancement by
declaring each closed vertex $x_j$ to be {\em green}
 if it is
correctly configured and $Z_j < q$.  We shall sometimes need
to consider the Poisson process with an extra vertex inserted
at $x \in B_n$, in which case the extra vertex 
 has values $Y_0$ and $Z_0$ associated with it. We shall refer to vertices that are either
red or green as being {\em coloured}.

Let $\partial B_n$ be the annulus $B_n \setminus B_{n - 0.5}$ and let
$A_n$ be the event that for the Poisson porcess $\Po_\lambda \cap B_n$,
 there is a path from a coloured vertex in
$B_{0.5}$ to a coloured vertex in $\partial B_n$ in $G
(\Po_\lambda,1)\cap B_n$  using only coloured
vertices (note that $A_n$ is based on a process completely
 inside $B_n$; we do not allow vertices outside of $B_n$ to
 affect possible enhancements inside $B_n$). 
For $x \in B_n$, let $A_n^x$ be defined the same way as $A_n$,
but in terms of the point process $(\Po_{\lambda} \cap B_n) \cup \{x\}$,
i.e. the Poisson process in $B_n$ with a point inserted at $x$. 

Let $\theta_n(p,q)$ be the probability that $A_n$ occurs, and define
$$\theta(p,q) \equiv \liminf _{n\rightarrow \infty} (\theta_n(p,q)).$$

The following proposition states that $\theta(p,q)$ is indeed the
percolation function associated to the enhanced model. From now on
we use `vertex' to refer to a point of the Poisson process and `point' to
refer to an arbitrary location in $\mathbb{R}^2$.

\begin{prop}\label{prop:enh}
There is a.s. an infinite connected component in $G(\Po_ \lambda,1)$ using only red and green vertices if and only if $\theta(p,q) > 0$.
\end{prop}
\noindent {\bf Proof of Proposition~\ref{prop:enh}.} For the if part let $A'_n$ be the event that there is a coloured path from $B_{0.5}$ to outside $B_{n-2}$, so $A_n$ is contained in $A'_n$. Let $\phi_n(p,q)$ be the probability of $A'_n$ occurring (which is monotone in $n$), and let $\phi(p,q)$ be the limit as $n$ goes to $\infty$. Therefore $\phi_n(p,q) \geq \theta_n(p,q)$ for all $n$ so $\phi(p,q) \geq \theta(p,q) > 0$, but $\phi(p,q)$ is just the probability of there being an infinite coloured component intersecting $B_{0.5}$ and it is well known that there is almost surely an infinite coloured component if 
$\phi(p,q) > 0$.

 For the only if part, if there is almost surely an infinite component then $\phi(p,q) > 0$.   Given $n \geq 6$, we build up the Poisson process on the whole of $B_{n - 3}$. If there are any closed vertices that are not definitely correctly or incorrectly configured, we build up the process in the rest of their $1$-neighbourhood, and this determines whether they are green or uncoloured. If any more closed vertices occur they cannot be correctly configured as they will be joined to a closed vertex. Therefore we have built up the process everywhere in a region $R $ with $B_{n-3} \subset R \subset B_{n-2}$,
 and all uncoloured vertices at this stage will remain uncoloured.
 Let $V$ be the set of coloured vertices that are joined by a coloured path to a coloured vertex in $B_{0.5}$ at this stage.

Next, we build out the process radially symmetrically from $B_{n-3}$ (apart from where the process has already been built up) until a vertex $v$ occurs that is connected to a vertex in $V$. Let $J$ be the
event that such a vertex $v$ occurs,
 so $J$ must occur for $A'_n$ to occur. 
Assuming $J$  occurs, set $r = |v|$, so
 $ r \in [n-3, n-1)$. Then
we can find points $a_1, a_2, \ldots, a_9$ on the line $0 v$ extended away from the origin such that $a_1$ is 
$r + 0.3$ from the origin, $a_2$ is $r+0.6$ from the origin and so on. Surround $a_1,\ldots, a_9$ with circles $D_1, \ldots, D_9$ of
radius $0.05$ around them. If there is at least one red vertex in
each one of these little circles that is contained in $B_n$
when the process continues to the whole
of $B_n$, and $v$ is also red then $A_n$ occurs. Therefore if $J$ occurs then the
conditional probability of $A_n$ occuring is at least $\gamma$,
where
\[
    \gamma = p(1 - \exp(- 0.0025 \lambda p \pi))^9,
\]
as this is the probability of getting at least one red vertex in
each little circle and $v$ being red. Therefore $\theta_n (p,q) \geq \gamma P[J] \geq
\gamma \phi (p,q)$ for all $n \geq 6$, so $\theta (p,q) \geq \gamma
\phi (p,q) > 0$. \hfill{$\Box$} \vspace{.5cm}

Our next lemma provides an analogue of the Margulis-Russo formula for
the enhanced continuum model. First, we need to introduce the notion of
pivotal vertices.

Given the configuration
$(\Po_{\lambda},\underline{Y},\underline{Z})$ and
inserting a vertex at $x$ we say that $x$ is $1$-$pivotal$ $in$
$B_n$ if putting $Y_0 = 0$ means that $A^x_n$ occurs but putting $Y_0
= 1$ means it does not. Notice that $x$ can either complete a path
(but it cannot do via being enhanced), or it could make another
closed vertex correctly configured which in turn would complete a
path. We say that $x$ is $2$-$pivotal$ $in$ $B_n$ if inserting a
vertex at $x$ and putting $Z_0 = 0$ means $A^x_n$ occurs but putting
$Z_0 = 1$ means it does not. That is, $Y_0 > p$ and adding a closed
vertex $v$ at $x$ means $v$ is correctly configured and enhancing it
to a green vertex means $A^x_n$ occurs but otherwise it does not.

For $i = 1,2$ let $E_{n,i}(x)$ be the event that $x$ is $i$-pivotal in $B_n$,
and set $P_{n,i}(x,p,q) := P[E_{n,i}(x)]$.

\begin{lemm} \label{prop:Russo}
For all $n > 0.5$ and $p \in (0,1)$ and $q\in (0,1)$ it is the case that
\bea
\frac{\partial \theta_n(p,q)}{\partial p}=\int^{}_{B_n} \lambda P_{n,1}(x,p,q) \, \mathrm{d}x
\label{eq1}
\eea
and
\bea
\frac{\partial \theta_n(p,q)}{\partial q} = \int^{}_{B_n} \lambda P_{n,2}(x,p,q) \, \mathrm{d}x.
\label{eq2}
\eea
\end{lemm}
\noindent {\bf Proof.}
Let ${\cal F}$  be the $\sigma$-algebra generated
by the locations but not the colours of the vertices of $\Po_\lambda \cap B_n$.
Let $N_1$ be the number of $1$-pivotal vertices.
Define ${\cal F}$-measurable random variables,
$X_{p,q}$ and 
$Y_{p,q}$ as follows; $X_{p,q}$ is
the conditional probability that $A_n$ occurs,
 and $Y_{p,q}$ is the conditional expectation of $N_1$, 
 given the configuration
of $\Po_\lambda$. 
By the standard version of the Margulis-Russo formula
for an increasing event defined on a finite collection
of Bernoulli variables (Russo~(1981), Lemma 3),
$$
\lim_{ h \to 0}
h^{-1} (X_{p+ h,q} - X_{p,q})
= 
 Y_{p,q} , ~~ a.s.
$$
Let  $M$ denote the total number of vertices of $\Po_\lambda$ in $B_n$.
By the standard coupling of Bernoulli variables, and
Boole's inequality,
$|X_{p+ h,q} - X_{p,q}| \leq |h| M$ almost surely, 
and since $M$ is integrable,
 dominated convergence yields
\begin{eqnarray}
\frac{\partial \theta_n(p,q)}{
\partial p}
=
\lim_{ h \to 0}
E[ h^{-1} (X_{p+ \delta,q} - X_{p,q})]
= 
E[
 Y_{p,q} ] = E[N_1],
\label{Mar29a}
\end{eqnarray}
and by a standard application of the Palm theory of Poisson processes
(see e.g. Penrose~(2003)),
the right  hand side of (\ref{Mar29a})
equals the right hand side
of (\ref{eq1}). The proof of (\ref{eq2}) is similar. 
\hfill{$\Box$} \\

The key step in proving Theorem $1$ is given by the following result.
\begin{lemm}\label{intermediate}
There is a continuous function $\delta:(0,1)^2 \to (0,\infty)$ such that
for all $n > 100$, $x \in B_n$ and $(p,q) \in (0,1)^2$, we have
\begin{equation}
P_{n,2}(x,p,q) \geq \delta(p,q)P_{n,1}(x,p,q).
\label{100419a}
\end{equation}
\end{lemm}

Before proving this, we give a result saying that we can
assume there are only red vertices inside an annulus of fixed
size. For $x \in  \R^2$, and $0 \leq \alpha < \beta$,   let 
$C_{\alpha}(x)$ be the closed circle (i.e., disc) of radius $\alpha$
centred at $x$, and let $A_{\alpha,\beta}(x)$ denote the annulus 
$C_{\beta}(x) \setminus  C_{\alpha}(x)$.
 Given $n$ and given $x \in B_n$, let 
$R_n(x,\alpha,\beta)$ be the event that all vertices in $A_{\alpha,\beta}(x) \cap B_n$
are red.

\begin{lemm}\label{intermed2}
Fix $\alpha >3$ and and $\beta > \alpha +3$.
There exists a continuous function $\delta_1:(0,1)^2 \to (0,\infty)$,
such that  
for all $(p,q) \in (0,1)^2$, all
  $n > \beta +3$ and all $x \in B_n$ with $|x|< \alpha -2 $ or 
$|x| > \beta +2$, we have
\[
P[E_{n,1}(x)  \cap R_n(x,\alpha, \beta) ] \geq \delta_1(p,q) P_{n,1}(x).
\]
\end{lemm}

\noindent {\bf Proof.}
 We shall consider a modified model, which is the same
as the enhanced model but with   enhancements suppressed
for all  those vertices lying in $A_{\alpha-1,\beta+1} := A_{\alpha-1,\beta+1}(x)$.
Let $E'_{n,1}(x) $ be the event that $x$ is 1-pivotal in
the modified model.

Returning to the original model,
 first create the Poisson process of intensity $\lambda$
in $B_n$.
Then for all the vertices  
in $B_n \setminus A_{\alpha-1,\beta +1}$,
 decide whether they are red or closed.
 Then, for  all those vertices in
 $B_n \cap A_{\alpha-1,\beta +1} $
with more than $4$ neighbours, or
with at least one closed neighbour outside $A_{\alpha-1,\beta +1}$,  decide
 whether they are red or closed. This decides whether or not they
 are coloured as these vertices cannot possibly become green because
 they are not correctly
 configured. We now can tell which of the closed vertices 
outside $A_{\alpha-1,\beta +1}$ are correctly configured, and we determine
which of these are green.

 This leaves a set $W$ of vertices inside $A_{\alpha-1,\beta +1}$ that have at most four neighbours. If we surround each vertex in $W$ by a circle of radius $0.5$ then we cannot have any point covered by more than $5$ of these circles as this means that there is a vertex in $W$ with at least $5$ neighbours. All
of these circles are contained in $C_{\beta+2}$,
 which has area $\pi(\beta+2)^2$. Therefore 
\begin{equation}
|W| \leq \frac{5 \pi (\beta+2)^2}{0.5^2 \pi} =  20(\beta + 2)^2.
\label{0513a}
\end{equation}
For $x$ to have any possibility of being $1$-pivotal,
at this stage there must be a
 set $W'$ contained in $W$ such that if every vertex in $W'$ is coloured and every vertex in $W \setminus W'$ is uncoloured then $x$ becomes $1$-pivotal.
In this case, with probability at least
$[p(1-p)]^{20(\beta +2)^2}$ we have every vertex in 
$W'$ red and every vertex in $W \setminus W'$ closed, which would imply 
 event $E'_{n,1}(x)$
 occurring. Therefore
 $P[E'_{n,1}(x)] \geq [p(1-p)]^{20(\beta+2)^2} P[E_{n,1}(x)]$.

Now we note that the occurrence or otherwise of $E'_{n,1}(x)$
is unaffected by the addition or removal of closed vertices
in $A_{\alpha, \beta}(x)$. This is because the suppression
of enhancements in $A_{\alpha-1,\beta+1}$ means that 
these added or removed vertices cannot be enhanced themselves,
and moreover any vertices they cause to be correctly
or incorrectly configured also cannot be enhanced.

Consider creating the marked Poisson process in $B_n$,
with each Poisson point (vertex) $x_i$  marked with 
the pair $(Y_i,Z_i)$, in two stages. First, add all marked vertices
in $B_n \setminus A_{\alpha,\beta}(x)$, and just the red vertices
in $B_n \cap A_{\alpha,\beta}(x)$. Secondly, add the closed vertices
in $B_n \cap A_{\alpha,\beta}(x)$.
 The vertices added at the second
stage have no bearing on the event $E'_{n,1}(x)$, so 
$E'_{n,1}(x)$ is independent of the event that no vertices 
at all are added in the second stage.
Hence, 
$$
P[E'_{n,1}(x) \cap R_n(x,\alpha,\beta)]
 \geq \exp(- (1-p) \lambda  \pi (\beta^2 - \alpha^2))
P[E'_{n,1}(x)], 
$$
with equality if $|x| \leq n- \beta$.

Finally, we use a similar argument to the initial argument in this proof.
Suppose
$E'_{n,1}(x) \cap R_n(x,\alpha,\beta)$ occurs. Then there exist at most 
$20(\beta +2)^2$ vertices in
 $A_{\beta,\beta +1}(x) \cup A_{\alpha -1,\alpha}(x)$
which are correctly configured for which the possibility of
enhancement has been suppressed. If we now allow these to be possibly enhanced,
there is a probability  of at least
$(1-q)^{20(\beta +2)^2}$ that none of them is enhanced, in which
case the set of coloured vertices is the same for the modified model
as for the un-modified model and therefore $E_{n,1}(x)$ occurs.  
Taking
 $$
\delta_1(p,q) = [p(1-p)(1-q)]^{20(\beta+2)^2} \exp(-(1-p) \lambda
 \pi (\beta^2 - \alpha^2 )),
$$
 we are done. \hfill{$\Box$} \\


\noindent {\bf Proof of Lemma \ref{intermediate}.}
Fix $p$ and $q$.  Also fix $n$ and  $x \in B_n$,  
and just write $P_{n,i}(x)$ for $P_{n,i}(x,p,q)$.
Define event $E_{n,1}(x)$
 as before,
so that $P_{n,1} (x) = P[E_{n,1}(x)]$. 
Also, for $0<r < s$ write
$C_{r}$ for the disc $C_r(x)$ and $A_{r,s}$ for the annulus $A_{r,s}(x)$.
For now we assume $30.5 < |x| < n-30.5$.
 We create the Poisson process of intensity $\lambda$
 everywhere on $B_n$ except inside
$C_{30}$, and decide which of these vertices are red.

Now we create the process of only the red vertices 
 in $A_{25,30}$ (a Poisson process of intensity $p \lambda$
in this region).
Assuming there will be no closed vertices in $A_{25,30}$,
we then know which of the closed  vertices outside $C_{30}$ are correctly
configured, and we determine which of these are green.

Having done all this, let  $V$ denote the set of current  vertices 
in $B_n \setminus C_{25}$ that are connected to $B_{0.5}$ at this stage
(by connected we mean connected via a coloured path), and
let $T$ denote the set of current vertices in $B_n \setminus C_{25}$ that
 are connected to $\partial B_{n}$.

Let $N(V)$ be the  $1$-neighbourhood of $V$ 
and  let $N(T)$  be the  $1$-neighbourhood of $T$. 
Recalling that $A \triangle B := (A \cup B) \setminus A \cap B$,
 we build up the red process inwards (i.e., towards $x$ from the boundary
of $C_{25}$)
on $C_{25} \cap (N(V) \triangle N(T))$
 until a red vertex $y$
occurs (if such a vertex occurs). Set  $r= |y-x|$.
Suppose $y \in N(V)$ (if instead $y \in N(T)$
we would reverse the roles of $V$ and $T$ in the sequel).
Then if $T \cap C_{r+0.05} \neq \emptyset$ we say that 
event $F$ has occurred and we let $z$ denote an
arbitrarily chosen vertex
of  $T \cap C_{r+0.05} $.
 Otherwise, we  build up the red process inwards on
 $C_{r} \cap N(T) \setminus N(V)$ until a red vertex $z$
occurs (if such a vertex occurs). 

Let $E_2$ be the event that (i)  such vertices $y$ and $z$ occur,
and (ii) the sets $V$ and $T$ are disjoint, and (iii)
 $|y-z| >1$, 
and (iv) there is no path from $y$ to $z$ through coloured
vertices in $B_n \setminus
C_{25}$ that are not in $V \cup T$.
If $E_{n,1}(x) \cap R_n(x,20,30)$ occurs, then $E_2$ must occur.

\begin{figure}[htbp]
\includegraphics[angle = 270, width = 14cm]{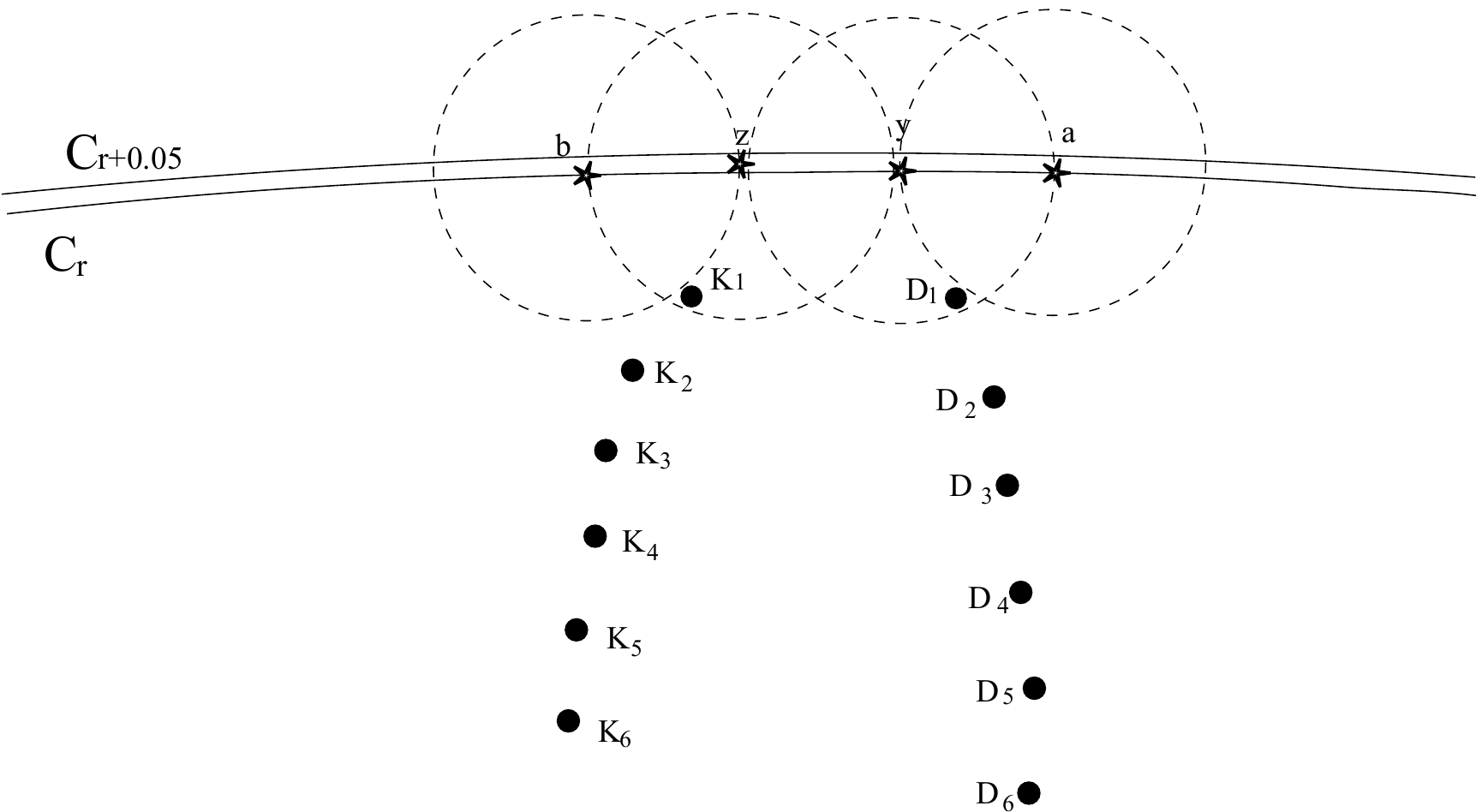}
\caption{Our convention in the diagrams is to indicate points with
lower case letters, and areas with upper case letters. The
dashed circles are of radius 1. Here the event $F$ occurs. 
} \label{fig:fig1}
\end{figure}

Now suppose $E_2 \cap F$ has occurred.  
Let $a$ be the point (again we use `point'
to refer to a point in $\mathbb{R}^2$) which is at distance
$r$ from $x$ and distance $1$ from $y$ on the opposite side of the
line $xy$ to the side $z$ is on
 (see Figure $\ref{fig:fig1}$). 
 Similarly let $b$ be the point lying at distance $1$ from $z$ and distance $r$ from $x$, on the opposite side of $xz$ to $y$.

Let $a_1$ be the point lying inside $C_r$ at distance 
$1.01$ from $a$ and $0.99$ from $y$, and
 let $D_1$ be the disc $C_{0.005}(a_1)$.
Let $b_1$ be the point at distance $1.01$ from $b$ and $0.99$ from
$z$, and let $K_1 := C_{0.005}(b_1)$.

Any red vertex in $D_1$ will be connected  to 
$y$ (and therefore to a path to $B_{0.5}$) but cannot be connected
to any coloured path to $\partial B_n$ as $a$ is the nearest place for
such a vertex to be, given $E_2 \cap F$ occurs.
Any red vertex in $K_1$ will be connected to $z$ (and therefore a path to
$\partial B_n$), but not a path to $B_{0.5}$. Also, any vertex in
$D_1$ will be at least $1.1$ away from any vertex in $K_1$. 

Now let $l$ be the line through $x$ such that $a_1$ and $b_1$ are on
different sides  of the line and at equal distance from the line.
We can pick points $a_2,a_3,\ldots,a_{30}$ such that $|a_i-a_{i-1}| \leq 0.9$
for $2 \leq i \leq 30$, and $\max(|a_{30}-x|, |a_{29}-x| ) \leq 0.9$,
 but  $|a_i -x| > 1.1$ for $i \leq 28$, and none of the
 ${a_i :i \geq 2}$ are within
$1$ of $C_r$ or within $0.51$ of $l$ or within $0.01$ of another
$a_j$. 

Do the same on the other side of $l$ with $b_2,b_3, \ldots , b_{30}$.
For $2 \leq i \leq 30$, define discs
 $D_i := C_{0.005}(a_i)$ and $K_i :=  C_{0.005}(b_i)$. 

Let $I$ be the event that there is at exactly one red vertex in each
of the circles $D_i$ and $K_i$, $1 \leq i \leq 30$,
and there are no more new vertices anywhere else in $C_{25}$, and no closed
vertices in $C_{30}\setminus C_{25}$. Then
\[
P[I|E_2 \cap F] \geq
(0.005^2 \pi \lambda p)
^{60} \exp(-900 \pi \lambda)
=: 
 \delta_2 .
\]
 If the events $E_2$,  $F$, $I$ occur and $Y_0
> p$ then $x$ is $2$-pivotal.

\begin{figure}[htbp]
\includegraphics[angle = 270, width = 12cm]{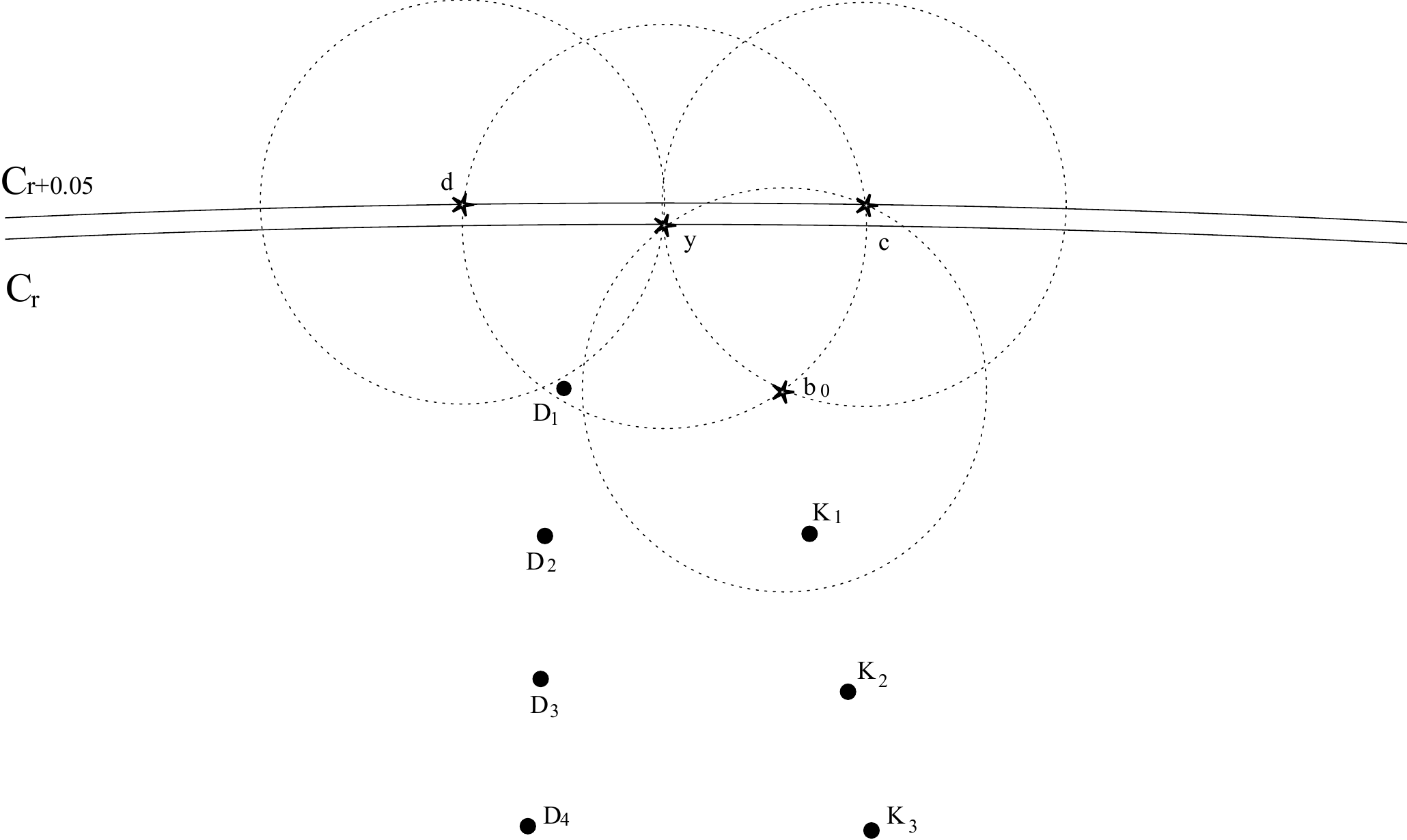}
\caption{The case where $F$ does not occur. Here $b_0$ is the `worst possible'
location for $z$.}
\label{fig:points1}
\end{figure}

Now we consider the case where $E_2$ occurs but $F$ does not, so $z$ is inside
$C_r$ and is connected to a vertex $z_1$ in $T$ that must be outside
$C_{r+0.05}$ because $T \cap C_{r+0.05} = \emptyset$
 (see Figure $\ref{fig:points1}$).

 Let $c$ be the point at distance $1$ from $y$ and $r+0.05$ from $x$, on
the same side of the line $xy$ as $z$ (assume without loss of
generality this is to the right of $y$). This is the closest $z_1$
can be. Let $b_0$ be the point inside $C_r$ at distance $1$ from $y$
and $1$ from $c$, so this is the furthest left that $z$ can be. Let
$d$ be the point at distance $r+0.05$ from $x$ and $1$ from $y$, on
the other side of $y$ to $c$. 
Let  $a_1$ be the point inside $C_r$ at
distance $1.01$ from $d$ and $0.99$ from $y$, and let
$D_1 := C_{0.005}(a_1)$. Then any vertex in $D_1$ is distant at
least $1.01$ from $b_0$, and therefore from $z$, as $z$ cannot be any
nearer than $b_0$. Also any vertex in $D_1$ will be at least $1.005$
from any other vertices in $T$, as $d$ is the nearest place such a
point can be. As before we can then have small discs $D_2,\ldots,D_{30}$
and $K_1,\ldots,K_{30}$ (of radius $0.005$)
 such that having one red vertex in each of these vertices ensures that
 $x$ is $2$-pivotal. 

Given $E_2 \setminus F$ occurs, the
probability of getting $1$ red vertex in each of the discs
$D_i$ and $K_i$ for $1 \leq i \leq 30$,
 and no other new vertices in $C_{25}$,
and no closed vertices in $C_{30} \setminus C_{25}$, is
at least $\delta_2$.
If this happens and also $Y_0>p$ then $x$ is $2$-pivotal.

So by Lemma \ref{intermed2},
 the probability that $x$ is $2$-pivotal satisfies
\beaa
    P_{n,2}(x)&  \geq &   \delta_2 (1-p) P[E_2  \cap F] +
      \delta_2 (1-p) P[E_2 \cap F^c] 
\\ 
& \geq &  \delta_2 (1-p)  P[ E_{n,1}(x) \cap R_n(x,20,30)]
\\ & \geq & \delta_1 \delta_2 (1-p)  P_{n,1}(x) .
\eeaa
 This proves the claim (\ref{100419a}) for the case with $30.5 < |x| < n-30.5$.

Now suppose $|x| \leq 30.5$.
 Create the Poisson process  in $B_n \setminus C_{40}$,
 and decide which of these vertices are red.  Then create  the
red process in  $A_{39,40}(x)$, and determine which
vertices in $B_n \setminus C_{40}$ are green, assuming
there are no closed vertices in $A_{39,40}(x)$.
Then build up the red process in $C_{39}$ inwards towards $x$
until a vertex $y$ occurs in the process which is connected to
$\partial B_n$. Let $H_1$ be the event that such a vertex $y$
appears at distance $r$ between $38$ and $39$ from $x$, so $H_1$ must
occur for $E_{n,1}(x) \cap R_n(x,20,40) $ to occur.

If $x$ is inside $B_{0.5}$ we can choose points $a_0$ and $a_1$ such
that they are both outside $B_{0.5}$, at distance between $0.8$ and
$0.9$ from $x$ and at distance between $0.1$ and $0.2$ from each
other. We can then choose $b_0$ and $b_1$ such that they are both
within $0.9$ of $x$, further than $1.5$ from $a_0$ and $a_1$ and
between $0.1$ and $0.2$ from each other. We can then choose points
$a_2, a_3, \ldots, a_{100}$ such that $|a_i -a_{i-1}| \leq 0.9$
for $2 \leq i \leq 100$, and $|a_{100}-y| \leq 0.9$, no two
$a_i$ are within $0.1$ of each other, and no $a_i$ is within $1.1$
of $x$, $b_0$ or $b_1$, or inside $B_{0.5}$ for $i \geq 2$. 

Define discs $D_i = C_{0.05}(a_i)$ and $K_j= C_{0.05}(b_j)$ 
 If there is at least one red vertex in each of these
discs and no vertices anywhere else in $C_r$, and $Y_0>p$,
 then $x$ is $2$-pivotal. If $x$ is outside $B_{0.5}$ we choose points in a
similar way but make sure $b_1$ connects with a path to $B_{0.5}$,
using little discs $K_2, K_3,\ldots,K_{50}$ which are again of radius
$0.05$ and are at least $1.1$ from the $a_i$. Therefore, 
setting
\[
\delta_3 := (1-p)(0.05^2 \pi \lambda p)^{152}\exp(-1600 \pi \lambda)
\]
and using
Lemma \ref{intermed2}, we have for some strictly positive continous
$\delta_4(p,q)$ that
$$
P_{n,2}(x) \geq \delta_3 P[H_1] \geq \delta_3 P[E_{n,1}(x)
  \cap
R_n(x,20,40)] \geq  \delta_3 \delta_4 P_{n,1}(x).
$$

Now suppose $|x| \geq n - 30.5$.
 In this case the proof is similar. Again, create
 the Poisson process in $B_n \setminus C_{40}$.
 Then create  the red process in $A_{39,40}(x)$ and
determine the colours of the vertices in $B_n \setminus C_{40}$,
assuming there are no closed vertices in $A_{39,40}(x)$.
Then  build the red process in $C_{39}  \cap  B_{n-0.5}$
  inwards towards $x$ until a vertex $y$ occurs that is connected to a
path of coloured vertices to $B_{0.5}$ but not to $\partial B_n$. Let
$H_2$ be the event that such a vertex $y$ occurs at distance $r$
between $38$ and $39$ from $x$, and that there is no current coloured
path from $B_{0.5}$ to $\partial B_n$, so $H_2$ must occur for 
$E_{n,1}(x) \cap R_n(x,20,40)$ to occur.
 Given this vertex $y$
 we can find discs
$D_1, D_2, \ldots ,D_{100}$ and $K_1, K_2, \ldots,K_{50}$ of radius $0.05$
as before such that having a red vertex in each of these discs
 but no other vertices in $C_r$ or $\partial B_n \cap C_{40}$,
and having $Y_0>p$,
ensures $x$ is $2$-pivotal. Therefore in this case
\[
P_{n,2}(x) \geq  \delta_3 P[H_2 ] \geq \delta_3 P[E_{n,1}(x) \cap
R_n(x,20,40)] \geq \delta_3  \delta_4
P_{n,1}(x).
\]

Take $\delta(p,q) := \min(\delta_1 \delta_2(1-p) ,\delta_3 \delta_4) $. By its
construction $\delta$ is strictly positive and continuous in $p$ and $q$,
and (\ref{100419a}) holds for all $x \in B_n$,
completing the proof of the lemma. \hfill{$\Box$} \vspace{.5cm}

The following proposition follows immediately from
Lemmas ~\ref{prop:Russo} and \ref{intermediate}.
\begin{prop} \label{prop:final}
There is a continuous function $\delta:(0,1)^2 \to (0,\infty)$
such that
for all $n \geq 100$ and $(p,q) \in (0,1)^2$, we have
\[
\frac{\partial \theta_n (p,q)}{\partial q} \geq \delta (p,q) \frac{\partial \theta_n (p,q)}{\partial p}.
\]
\end{prop}

\noindent
\textbf{Proof of Theorem~\ref{bool}.}
Set $p^* = p_c^{\rm site} $
and $q^* = (1/8)(p^*)^2$. 
Then using Proposition~\ref{prop:final} and looking at a small
 box around $(p^*,q^*)$, we can find 
$\epsilon \in (0,\min(p^*/2,1-p^*))$ 
 and $ \kappa \in (0,q^*)$ 
such that for all $n > 100$ we have
\[
\theta_n(p^* + \epsilon,q^* - \kappa) \leq
 \theta_n(p^*  - \epsilon,q^* + \kappa).
\]
Taking the limit inferior as $n \rightarrow \infty$,
since
 $\theta$ is monotone in $q$ we get
\[
0 < \theta(p^* + \epsilon,0 ) \leq
\theta(p^*  + \epsilon,q^*- \kappa)
 \leq \theta(p^* - \epsilon,q^* + \kappa).
\]
Now set $p = p^* - \epsilon $. Then   $q^* + \kappa \leq p^2$,
so that $\theta(p,p^2)>0$, and
 by Proposition
\ref{prop:enh}, the enhanced model with parameters $(p,p^2)$ percolates,
i.e. has an infinite coloured component, almost surely.

We finish the proof with a coupling argument along the lines of
Grimmett and Stacey~(1998). Let E be the set of edges and $V$ be the
set of vertices of $\CC$ (the infinite component). Let $(X_e: e \in
E)$ and $(Z_v: v \in V)$ be collections of independent Bernoulli
random variables with mean $p$. From these we construct a new
collection $(Y_v: v \in V)$ which constitutes a (red) site percolation
process on $\CC$, as follows. Let
 $e_0, e_1, ...$ be an enumeration of the edges of $\CC$ and
 $v_0, v_1, ...$ an enumeration of the vertices. Suppose
at some point we have defined $(Y_v: v \in W)$ for some subset $W$
of $V$. Let $\cal Y$ be the set of vertices not in $W$ which are
adjacent to some currently active vertex (i.e. a vertex $u \in W$
with $Y_u = 1$). If $\cal Y = \emptyset$ then let $y$ be the first
vertex not in $W$ and set $Y_y = Z_y$ and add $y$ to $W$. If $\cal Y \neq \emptyset$,
we let $y$ be the first vertex in $\cal Y$ and let $y'$ be the first
currently active vertex adjacent to it, then set $Y_y = X_{yy'}$ and add $y$ to $W$.
Repeating this process builds up the entire red site percolation process,
if it does not percolate, or a percolating subset of the red site percolation
process if it does percolate.
In the latter case, the bond process $\{X_e\}$ also percolates.

Now suppose the red site process does not percolate.
For any correctly configured vertex $x$ with $v,w,y,z$ as
in Figure \ref{fig:enh1}, $x$ itself is not red. Therefore at most one edge to $x$
has been examined, so we can can find a first unexamined edge (in
the enumeration) to $v$ or $w$, and then to $y$ or $z$. We
then declare $x$ to be green only if both of these edges are open,
which happens with probability $p^2$.
This completes the enhanced site process with
 $q = p^2$ and every component
in this is contained in a component for the bond process $\{X_e\}$.

Therefore, since the enhanced $(p,p^2)$ site process percolates almost
surely,
so does the bond process, 
so $p_c^{\rm bond} \leq p < p_c^{\rm site}$. 
\hfill{$\Box$}
\vspace{.5cm}

\section{RCM: the key lemma}
\label{seclemma}
\eqco

This section is devoted to stating and proving Lemma
 \ref{Hprop} below, which
is the key step in subsequently 
  proving Theorems \ref{RCM}
 and \ref{th:squash}.
We consider the RCM with connection function 
 $f:[0,\infty) \rightarrow [0,1]$. Throughout this section we
assume that $f$ is nonincreasing and, moreover, that 
\begin{equation}
 \sup \{a : f(a) > 0\} =1.
\label{rangeone}
\end{equation}
Fix $x \in \R^2$ and (as in the preceding section)
for $r <s$ let $C_r$ denote the disc of radius $r$ centred at $x$ and
let $A_{r,s}$ denote the annulus $C_s \setminus C_r$.

We consider the RCM on a Poisson process 
 in $C_{29}$, under certain {\em boundary conditions},
represented by three finite disjoint 
sets $V, T$ and $S$ in $\R^2 \setminus C_{29} $, together with a
 collection ${\cal E}$ of edges amongst the vertices (i.e., elements) of
$S$. We write ${\bf S}$ for the graph  $(S,{\cal E})$
(a subgraph of the complete graph on vertex set $S$). 
 We refer to the triple $(V,T,{\bf S})$
as a {\em boundary condition}.

In terms of generalising the proof of Theorem \ref{bool} to the RCM,
the set $V$ (respectively $T$)
 represents the set of coloured vertices in $B_n \setminus C_{29}$
 that are connected by a coloured path to $B_{0.5}$ (respectively, to
$\partial B_n$),
before the vertices inside $C_{29}$ have been added.
The set $S$ represents the remaining coloured vertices
$B_n \setminus C_{29}$, and  ${\cal E}$ represents the
set of edges between these vertices. However, this description
is only for motivation, and the present section
is self-contained; in particular, no colouring of vertices
takes place in this section.

For $\mu >0$ and $0 \leq r < s$, let $\Po_{\mu,r,s}$
denote a homogeneous Poisson process of intensity $\mu$ in $A_{r,s}$.
Given $(V,T,\bS)$ as described above, 
for $0 \leq r < 29$ the RCM on $\Po_{\mu,r,29}$ 
with boundary condition $(V,T,\bS)$ is  obtained as follows:
 we  connect each
 pair of vertices $x,y$ with $x,y  \in \Po_{\mu,r,29}$
or $x \in \Po_{\mu,r,29}$ and $y \in V \cup T \cup S$,
by an undirected edge
  with probability $f(|x-y|)$, 
independently of other pairs.
For $x \in \Po_{\mu,r,29}$ we then say $x$ is {\em path-connected}
to $T$ (respectively, to $V$) if there is a path
 from $x$ to $T$ (respectively, $V$) using the  edges created. 
If also $y \in \Po_{\mu,r,29}$ then we say $x$ is path-connected
to $y$ if there is a path from $x$ to $y$, using the edges
created along with the edges of ${\cal E}$.

 Let $V_r$, respectively $T_r$ be the set of vertices
of $\Po_{\mu,r,29}$
that are path-connected to $V$, respectively
$T$. 
Let $S_{r}$ be the set $\Po_{\mu,r,29} \setminus (V_r \cup T_r)$.
Define the event
\begin{equation}
H(V,T,\bS) : = \{V_{20} \cap C_{21} \neq \emptyset\} \cap 
\{T_{20} \cap C_{21} \neq \emptyset\} \cap 
\{ V_{20} \cap T_{20} = \emptyset \}.
\label{Hdef}
\end{equation}
Let $H'(V,T,\bS)$ be the intersection of $H(V,T,\bS)$
with the event that there exists $v^* \in C_{20.1}$
and $t^* \in  C_{20.1}$ such that $|v^*-t^*| >1.5$ and
$V_{20} \cap C_{20.5} = \{v^*\}$ and 
$T_{20} \cap C_{20.5} = \{t^*\}$, and $S_{20} \cap C_{20.5} = \emptyset$.
We can now state the main result of this section.

\begin{lemm}
\label{Hprop}
Suppose $f$ is nonincreasing and 
(\ref{rangeone}) holds.
Then there exists a continuous  function
 $\epsilon: (0,\infty) \to (0,\infty)$ such that
 for any $\mu \in (0,\infty)$,
and any boundary condition $(V,T,\bS)$
 we have 
\begin{eqnarray}
P[H'(V,T,\bS)] \geq \epsilon (\mu) P[H(V,T,\bS)].
\label{Hpropeq}
\end{eqnarray}
\end{lemm}

We shall need several further lemmas to prove
 Lemma $\ref{Hprop}$.
In these  arguments, we often need to build up the Poisson
 process $\Po_\mu$ in certain regions via  a ``scanning process", as  described in Meester, Penrose and Sarkar
 (1997) which gives a rigorous proof that it does indeed build
 up the Poisson process. For any set of vertices $U$ and
 any point $z \in \R^2$ let $p(z,U)$ be the probability that a vertex at
 $z$ is joined to at
least one of the vertices in $U$. So $1 - p(z,U) = \Pi_{u \in U}(1 -
f(|z - u|))$.

We shall consider the  process  $\Po_{\mu,24,25}$ as
 the union of two independent half-intensity processes 
$\Po_{\mu/2,24,25} $ and $\Po'_{\mu/2,24,25}$.
Let $E_1$ be the event that $\Po_{\mu /2,24,25}$
has precisely two elements, and one of these is connected to $V_{25}$
 while the other is connected to $T_{25}$,
and $V$ is not path-connected to $T$ through $\Po_{\mu,25,29}
\cup  \Po_{\mu/2,24,25} \cup S$.  

\begin{lemm}
\label{lemhalf1} For all boundary conditions $(V,T,\bS)$,
it is the case that $P[E_1] \geq 0.25 \exp(-25 \pi \mu) P[H(V,T,\bS)]$.
\end{lemm}
{\bf Proof.}
Create the  process $\Po_{\mu,25,29}$ 
and define the sets $V_{25}$,  $T_{25}$ and $S_{25}$ as described earlier.
Then build up an inhomogenous process in from the edge of $C_{25}$
 (i.e. starting at distance $25$ from $x$ and working radially symmetrically
 inwards) with intensity $\mu h_1(\cdot)$ where $h_1(v) = p(v,V_{25})
(1-p(v, T_{25}))$,
 until a vertex $y$ occurs.
Then add edges from $y$ to $V_{25}$ conditional on
 there being at least one such edge.
Add  edges independently from $y$ to vertices in $S_{25}$
 in the usual way.
Do not add any edges from $y$ to $T_{25}$.

Now build up another inhomogenous process in from the edge of $C_{25}$
 with intensity $\mu h_2(\cdot)$,  where 
$h_2(v) = p(v, T_{25}) (1-p(v,V_{25} ))$, until a vertex $z$ occurs.  
Add edges from $z$ conditional on there being at
 least one  edge from $z$ to $T_{25}$ and no edge from $z$ to $V_{25}$.

 Let $E'_1$ be the event that we get such vertices $y$ and $z$
 and $y$ is not connected to  $z$ through $S_{25}$.
Then $E'_1$ must occur for the event $H(V,T,\bS)$ to occur.

Let $E''_1$ be the event that $E'_1$ occurs with both   $y$ and $ z$ coming
 from the first half intensity
process $\Po_{\mu/2,24,25}$ (rather than from $\Po'_{\mu/2,24,25}$).
Then $P[E''_1 | E'_1] = 0.25$.
 Given $E''_1$ occurs, for $E_1$
to occur we need only that there be
 no further  vertices of
$\Po_{\mu/2,24,25}$ besides $y$ and $z$, and the conditional
probability of this is at least
$ \exp (- 49 \pi \mu/2)$.
Combining these probability estimates gives the result. 
\hfill{$\Box$} \\

Let $\rho := \inf \{a>0 : f(a) < 1\}$, i.e. the radius of certain
connection (this could be zero).
We shall prove Lemma \ref{Hprop} separately for the
two cases $\rho  < \frac{1}{ \sqrt{2}} - 0.01$ and
$\rho  \geq \frac{1}{ \sqrt{2}} - 0.01$ (see Lemmas \ref{Hcase1} and \ref{Hcase2} below).

\begin{figure}[htbp]
\includegraphics[angle = 0, width = 16cm]{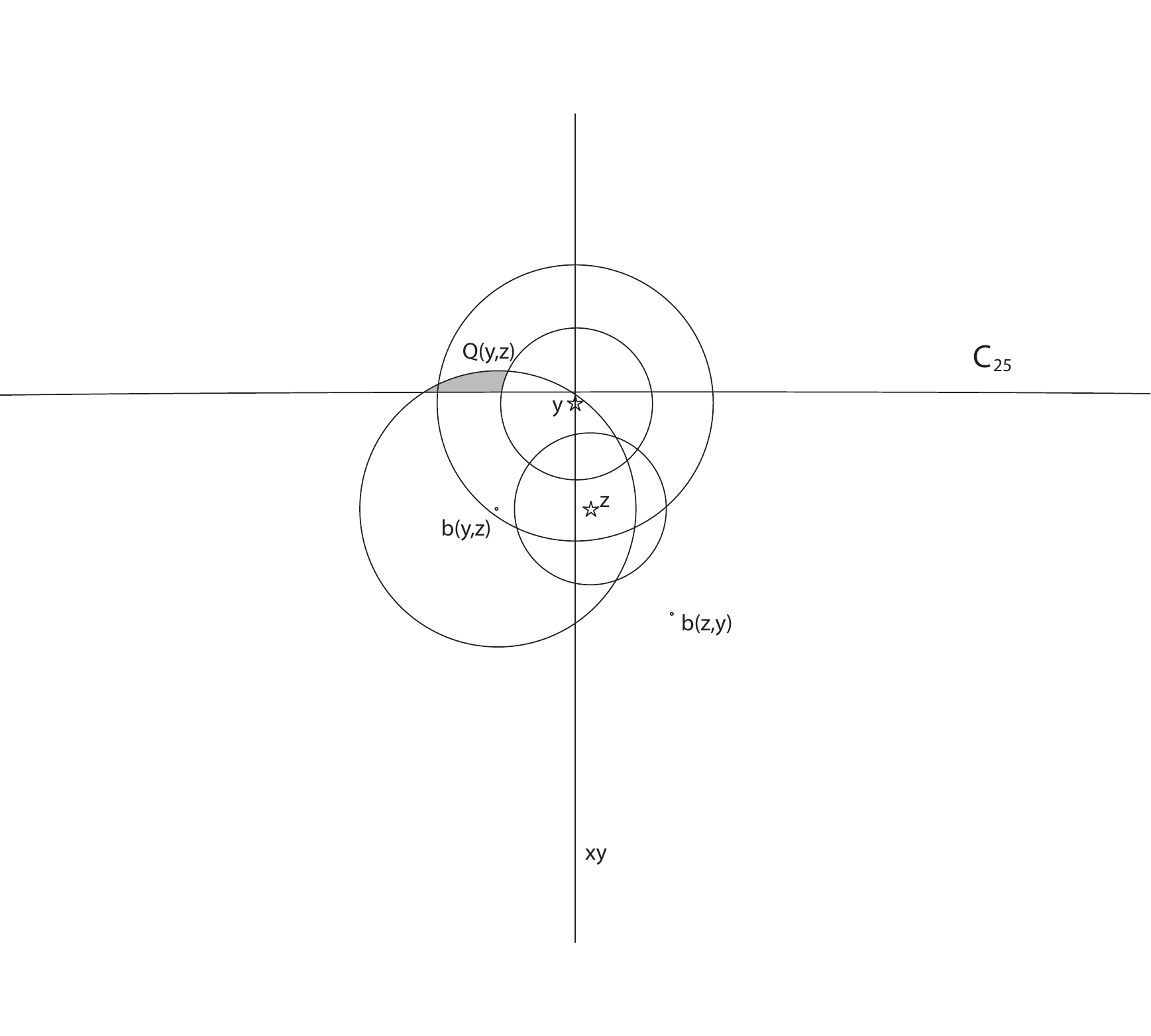}
 \caption{Here is a diagram showing the region $Q(y,z)$ (in this case
$Q(z,y)$ is empty). The smaller circles are of
 radius $\rho$ and the larger ones are of radius 
$1.0001$
} \label{fig:points2}
\end{figure}

Suppose for now that $\rho  < \frac{1}{ \sqrt{2}} - 0.01$.
Given $y,z \in A_{24,25}$ with $x,y,z$ not collinear, 
let $ b(y,z)$ be the point at distance 
$0.999$ from $y$, at distance
 $\rho + 0.01$ from $xy$ and on the opposite side of the line $xy$ to $z$ 
(see Figure \ref{fig:points2}).
Let  $ b(z,y)$ be defined similarly.
Define the region 
$$
Q(y,z) := C_{1.0001}(b(y,z)) \setminus ( C_{25} \cup C_{\rho}(y) )
$$
and define $Q(z,y)$ similarly
 (see Figure $\ref{fig:points2}$, where $Q(z,y)$ is empty). 
The regions $Q(y,z)$ and $Q(z,y)$, if non-empty,  each
 have diameter less than $0.9$
 due to $\rho$ being less than $\frac{1}{ \sqrt{2}} - 0.01$.

Given $y$ and $z$
define $T^{y,z}_{25}$ and $V^{y,z}_{25}$ in the
 same manner as $T_{25}$ and $V_{25}$, respectively,
but using the point process $\Po_{\mu,25,29} \cup \{y,z\}$
 instead of $\Po_{\mu,25,29}$.

Suppose $E_1$ occurs, and let $y, z$ be  the vertices of $\Po_{\mu/2,24,25}$,
 with $y$ path-connected to $V$ and $z$ path-connected to $T$. 
Let  $E_2$ be the event that there are no more than two vertices of 
$T^{y,z}_{25}$ in $Q(y,z)$ and no more than two vertices of 
$V^{y,z}_{25}$ in  $Q(z,y)$,
and no vertices of $\Po_{\mu,25,29}$ at all, other than those of $T^{y,z}_{25}$ and $V^{y,z}_{25}$.

\begin{lemm}
\label{Q}
Suppose $\rho < \frac{1}{\sqrt{2}} - 0.05$. Then
\[
P[E_2 | E_1 ] \geq 
f(0.9)^2 \exp(-29^2 \pi \mu)
= : \epsilon_1(\mu).   
\]
\end{lemm}

\noindent {\bf Proof.}
The idea here is to condition on what happens inside the
annulus
 $A_{24,25}$.
The probability $P[E_1]$ is the product of the probability that
 there are  exactly two vertices in $\Po_{\mu/2,24,25}$,
 and the probability that for two uniformly 
distributed vertices in $A_{24,25}$, they are joined one of them to 
$T_{25}$  but not $V_{25}$ and the other to $V_{25}$ but not $T_{25}$. 
Given $y$ and $z$ in $A_{24,25}$,  let $I_{y,z}$ be the event
(defined in terms of the Poisson process  $\Po_{\mu,25,29}$ and
associated edges)
that  $y \in V_{25}^{y,z} \setminus T_{25}^{y,z}$ and
  $z \in T_{25}^{y,z} \setminus V_{25}^{y,z}$, and
 let $p(y,z) = P[I_{y,z}]$ (this also depends on $V,T$ and $\bS$).
Then
$$
P[E_1] = \exp(-49  \pi \mu/2) (\mu/2)^2
 \int_{A_{24,25}}  \int_{A_{24,25}} p(y,z) dy dz. 
$$
Similarly,
$$
P[ E_1 \cap E_2]
 = \exp(-49   \pi \mu/2) (\mu/2)^2
 \int_{A_{24,25}}  \int_{A_{24,25}} p'(y,z) dy dz, 
$$
where $p'(y,z) = P[ I'_{y,z}]$ and $I'_{y,z}$ is the event that 
$I_{y,z}$ occurs and also there are at most two  vertices 
of $T^{y,z}_{25}$ in $Q(y,z)$,
and at most two  vertices of $V_{25}^{y,z}$ in $Q(z,y)$,
and all vertices in $\Po_{\mu,25,29}$
are in $V_{25,29}^{y,z} \cup T_{25,29}^{y,z}$. 
 Therefore we just need to show that $p'(y,z) \geq \epsilon_1 p(y,z)$ for
Lebesgue-almost all 
$y,z$ in $A_{24,25}$, and for all possible configurations where $E_1$ occurs. 
We do this in stages. \\

{\bf Stage 1.}
Fix $y$ and $z$.
Let $V_0 = V \cup \{y\}$ and $T_0 = T \cup \{z\}$.
We now 
{\em exhaustively create} the set of vertices in
 $A_{25,29} \setminus Q(z,y)$ that are path-connected
 to $V_0$ but not to $T_0$, by
which we mean the following sequence of steps. 
First create a process of intensity 
$\mu p(\cdot, V_0 )(1-p(\cdot,T_0))$ in 
$A_{25,29} \setminus Q(z,y)$.
 Add edges from the new vertices to $V_0 \cup S$  conditional
on  having at least one edge from each new vertex  to  $V_0$
and no edges from the new vertices to $T_0$.
Let $V_1$ be the set of vertices outside $V_0$ that are 
now path-connected to $V_0$
(i.e. the newly added vertices and any vertices of $S$ that are
 path-connected to them). 
Next, create a process of intensity 
$\mu p(\cdot,V_1)(1- p(\cdot, V_0)) (1 - p(\cdot,T_0))$ in 
$A_{25,29} \setminus Q(z,y)$, and add edges to these points
conditional on having at least one edge from each new
point to $V_1$ but no edge to $V_0$ or $T_0$. Let $V_2$ be
the set of points now path-connected to $V_0$ that were
not in $V_0 \cup V_1$. Next create a process in $A_{25,29}
\setminus Q(z,y)$ of
intensity $\mu p(\cdot,V_2)(1- p(\cdot,V_0 \cup V_1))(1-p(\cdot,T_0))$.  

Continue in this way, at each stage adding those
 vertices in $A_{25,29} \setminus Q(z,y)$
that are connected to the latest $V_i$ but not
to earlier sets $V_{i-1},\ldots,V_0$ or to $T_0$.
At some stage this procedure must terminate (i.e. the new
Poisson process has no points).
This completes the exhaustive creation of points that are path-connected
to $V_0$ but not $T_0$. 

Now let  $V'$ be the union of $V$ with all vertices path-connected to
 $V$ at this stage,
and let $U_y$ be the union of $\{y\}$ with the
 set of vertices path-connected to $y$ at this stage.

{\bf Stage 2.}
Next, we exhaustively create
 (in a similar manner to the above)
 the set of vertices in $A_{25,29} \setminus Q(y,z)$
that are path-connected to $T_0 $ but not to $V' \cup U_y$.
Then let $T'$ be the union of $T$ with all vertices path-connected  to $T$ 
at this stage, and let $U_z$ be the union of $\{z\}$ with
the set of all vertices path-connected 
to $z$ at this stage.

{\bf Stage 3.}
Suppose next that $z \notin T'$. Otherwise, go on to Stage 4 below.
 Then,
since we have exhaustively created the vertices connected
to $T' \cup U_z   $ outside $Q(y,z)$,
 for $I_{y,z}$
to occur there must be a vertex in $Q(y,z)$ connected to $T'$
and a vertex (possibly the same one)  in $Q(y,z)$ connected to $U_z$.
Build up the process in $Q(y,z)$ towards $x$ with intensity 
\[
\mu p(\cdot,U_z)p(\cdot,T') [1-p(\cdot,V' \cup U_y)]
\]
until we get a vertex $u$ (if any). If such a vertex occurs then we
 add edges from $u$ to $U_z$ and to $T'$ conditional on there being
 at least one of each type, and add no edges from $u$ to $V' \cup U_y$.
We then let
 $T'' := T' \cup U_z \cup \{u\}$,
 and go to Stage 4 below.

If $u$ does not occur, build up two more processes in
 $Q(y,z)$, with intensities
\[
\mu [1-p(\cdot,U_z)]p(\cdot,T') [1-p(\cdot,V' \cup U_y)]
\]
and
\[
\mu 
p(\cdot,U_z)[1-p(\cdot,T')] [1-p(\cdot,V' \cup U_y)]
\]
until we get vertices $u_1$ and $u_2$ respectively. If
 we get such vertices then $u_1$ will be joined to $T'$ and $u_2$ will be joined
 to $U_z$. Also, $u_1$ will be joined to $u_2$ with probability at least
 $f(0.9)$.
Assume this happens (so now $z$ is path-connected to $T$), and
 let $T'':= T' \cup U_z \cup \{u_1,u_2\} $
and go to Stage 4.
If we do not get $u_1$ and $u_2$, then $I_{y,z}$ cannot occur.

{\bf Stage 4.}
Suppose now that $y \notin V'$. Otherwise, go on to Stage 5 below.
Build up the process in $Q(z,y)$ towards $x$ with intensity 
\[
\mu p(\cdot,U_y)p(\cdot,V') [1-p(\cdot,T'' \cup U_z)]
\]
until we get a vertex $w$. If such a vertex occurs, then add edges
 from $w$ to $U_y$ and to $V'$ conditional on there being at least one of
 each type, add none to $T' \cup U_z$.
We now have a path from $y$ to $V$ and go to Stage 5 below.

If $w$ does not occur, build up two more processes in $Q(z,y)$, with
 intensities
\[
\mu 
[1-p(v,U_y)]p(v,V') 
[1-p(v,T'' \cup U_z)]
\]
and
\[
\mu
p(v,U_y)[1-p(v,V')] 
 [1-p(v,T'' \cup U_z)]
\]
until we get vertices $w_1$ and $w_2$ respectively. If
 we get such vertices, then $w_1$ will be
 joined to $V'$ and $w_2$ will be joined
 to $U_y$. Also $w_1$ will be joined to $w_2$ with probability
 at least $f(0.9)$.
Assume this happens (so  then we have a path from $y$ to $V$),
and go to Stage 5. If $w_1$ and $w_2$ do not occur, then 
$I_{y,z}$ cannot occur.

{\bf Stage 5.}
 By now we have $y$ connected (by a path) to $V$ and $z$ connected
to $T$, and $V$ not connected to $T$.
Now sample the rest of $\Po_{\mu,25,29}$.
Then as long as no more vertices occur when we do this
(an event with probability at least $\exp(-29^2 \pi \mu)$), 
event $I_{y,z}$
occurs. Therefore, 
we have shown that
$p'(y,z) \geq
 \epsilon_1 p(y,z)$, as required.
 \hfill $\Box$
\\

\begin{lemm}
\label{Hcase1}
Suppose that $f$ is nonincreasing and (\ref{rangeone}) holds, and
 that $\rho < \frac{1}{\sqrt{2}} - 0.05$. Then the conclusion
of Lemma \ref{Hprop} holds. 
\end{lemm}
{\bf Proof.}
Suppose $E_1 \cap E_2$ occurs, and let $y$ and $z$ be as in the
definition of $E_1$ (i.e. the points in $\Po_{\mu/2,24,25}$
that are path-connected to $V$ and to $T$ respectively).
Let $b_1 = b(y,z)$ and $a_1 = b_1(z,y)$.
 Define discs
 $D_1 := C_{0.0001}(b_1)$ and $K_1:= C_{0.0001}(a_1) $.
Then
\begin{eqnarray}
\min ( \dist(D_1,z) , \dist (K_1,y), \dist(D_1,K_1) )
 \geq \rho+0.005 ;
\label{0512a}
 \\ \min(
\dist(D_1,\R^2 \setminus C_{25}), \dist(K_1,\R^2 \setminus C_{25}) )
 \geq  \max(\rho + 0.005, 0.6),
\label{0625a}
\end{eqnarray}
and
 for any $b' \in D_1$ and $a' \in K_1$ we have
$\max(|b'-y|,|a'-z|) \leq 0.9991$.

Next take further discs $D_i = C_{0.0001}(b_i)$ and $K_i= C_{0.0001}(a_i)  $,
 for $2 \leq i \leq 7$, such
that each of these discs is contained in  
$A_{20,24}$, and discs $D_1,K_1,\ldots,$ $D_7,K_7$ are disjoint, and
\begin{eqnarray}
  |b_i-b_{i-1}|= |a_i-a_{i-1}| =  0.999,  ~~~  2 \leq i \leq 7;
\nonumber
\\
 \min(\dist(D_2,K_1), \dist(K_2,D_1)) \geq \rho +0.005 ;
\nonumber
\\
 \dist(D_i,K_j) \geq 1.1,  ~~~   1 \leq i,j \leq 7 , ~(i,j) \notin
\{(1,1),(1,2),(2,1)\};
\nonumber
\\
  \min(|b_i-x|,|a_i-x|) \geq 20.6, ~~~ 2 \leq i \leq 6; 
\nonumber
\end{eqnarray}
and $ |b_7-x| = |a_7-x|= 20.05$
and $|b_7-a_7| \geq 1.5$.  

Now create the Poisson process $\Po'_{\mu/2,24,25} \cup \Po_{\mu,20,24}$.
Let $E_3$ be the event that we get exactly one new vertex in each
of $D_i$ and $K_i$ (denoted $y_i$ and $z_i$ respectively) for $1 \leq i \leq 7$,
and
 no other new vertices.
Then
\begin{equation}
P[E_3 | E_1\cap  E_2] \geq ((0.0001)^2\pi \mu/2)^{14}
 \exp (-25^2 \pi \mu) =: \epsilon_2.
\label{0514a}
\end{equation}
Now, assuming $E_1 \cap E_2 \cap E_3$ occurs, decide which   edges 
occur involving
the new vertices. The probability that we get edges forming
the  paths
$(y,y_1,y_2,\ldots,y_7)$ and $(z,z_1,\ldots,z_7)$
is at least $f(0.9991)^{14}$.

By  (\ref{0512a}), the probability that $y_1$ is not joined
 to $z $,  $z_1$ or  $z_2$ is at least 
$[1-f(\rho+0.005)]^{3}$. Also by \eq{0625a},
 the probability that $y_1$ is joined to no vertices
 of $T_{25}^{y,z} \cap A_{25,29}$ is at least
$[1-f(\rho+0.005)]^{2}$, because
at most 2  such vertices lie in  $Q(y,z)$  
since event $E_2$ is assumed to occur,
and no such vertices  lie within  $\rho$ of $y$ since
event $E_1$ is assumed to occur, and all other such vertices
are more than unit distance from $y_1$.

Similarly, $z_1$ is not connected to  $y$ or $y_2$ or 
any vertex in $T_{25}^{y,z} \cap A_{25,29}$
with probability at least 
$[1-f(\rho+0.005)]^{4}$ given $E_1 \cap E_2$. 

If $y_2$ is not connected to $z_1$ and $z_2$ is not connected
to $y_1$, then
for $2 \leq i \leq 7$,
none of the vertices $y_i $ can be connected to any of
the vertices $z_j$ or to $T_{25}^{y,z} \cap A_{25,29}$,
and none of the vertices $z_i$ can be connected to any of
the vertices $y_j$ or to $V_{25}^{y,z} \cap A_{25,29}$.
Therefore,  we arrive at
$$
P[H'(T,V,\bS)|E_1\cap E_2 \cap E_3] \geq 
f(0.9991)^{14}
[1-f(\rho+0.005)]^{9}:= \epsilon_3.
$$ 
Hence, by (\ref{0514a}) and Lemmas \ref{lemhalf1} and \ref{Q},
taking
$
\epsilon = 0.25 \exp(-25 \pi \mu)
\epsilon_1\epsilon_2\epsilon_3,
$
 we have the desired result (\ref{Hpropeq})
 for $\rho \leq \frac{1}{\sqrt{2}}-0.05$.
\hfill $\Box$ \\

 Now, to complete the proof of Lemma \ref{Hprop} we
 look at the case where $\rho
\geq \frac{1}{\sqrt{2}}-0.05$.  We create the
 process
$\Po_{\mu,25,29}$
and define $V_{25}, T_{25} $ and $S_{25}$ as before.
Let $E_4$ be the event that $V_{25}$ and $T_{25}$ are disjoint.
This must occur for $E_1$ to occur.

We then add the half intensity process 
$\Po_{\mu/2,24,25}$.
Let $F_V$ be the event that $E_4$ occurs and there is just one vertex $y$ of 
$\Po_{\mu/2,24,25}$, and it is connected to $V_{25}$
but not $T_{25}$, and $T_{25}$ includes a vertex in $A_{|y-x|,|y-x|+.05}$.  
Similarly,
let $F_T$ bve the event that $E_4$ occurs and there is just one vertex $y$ of 
$\Po_{\mu/2,24,25}$, and it is connected to $T_{25}$
but not $V_{25}$, and $V_{25}$ includes a vertex in $A_{|y-x|,|y-x|+.05}$  

Let $G_V$ be the event that $E_4$ occurs and there are just
 two vertices $y,z$ of
$\Po_{\mu/2,24,25}$, and $y$ is connected to $V_{25}$
but not $T_{25}$
 and $z$ is connected to $T_{25}$ but not $V_{25}$, and  $|y-x| > |z-x|$ and
$T_{25} \cap A_{|y-x|,|y-x|+.05} = \emptyset$. Similarly  
let $G_T$ be the event that $E_4$ occurs and there are two vertices $y,z$ of
$\Po_{\mu/2,24,25}$, and $y$ is connected to $T_{25}$
but not $V_{25}$
 and $z$ is connected to $V_{25}$ but not $T_{25}$, and  $|y-x| > |z-x|$ and
$V_{25} \cap A_{|y-x|,|y-x|+.05} = \emptyset$.

\begin{lemm}
\label{lem0507}
Let $\epsilon_4(\mu) :=0.25\exp(-25 \pi \mu)$.
Then for any boundary conditions $(V,T,\bS)$ we have 
\begin{equation}
P[H(V,T,\bS)] \leq 
\epsilon_4^{-1} ( P[F_V] + P[ F_T] + P[G_V] + P[G_T]).
\label{0507a}
\end{equation} 
\end{lemm}
{\bf Proof.}
After creating the process $\Po_{\mu,25,29}$,
we build  a process of intensity 
$$
\mu (p(\cdot,V_{25})(1- p(\cdot,T_{25}))
 + p(\cdot,T_{25})(1- p(\cdot, V_{25})))
$$
 inwards into $C_{25}$, until we get a vertex $y \in A_{24,25}$. 
Let
 $E'$ be the event that such a vertex occurs. 
Event $E'$ must occur for $H(V,T,\bS)$ to occur.

If $E'$ occurs,
add edges from $y$ to $V_{25} \cup T_{25} \cup S_{25}$, conditional
on there being at least one edge from $y$ to $V_{25} \cup T_{25}$ but
there not being edges from $y$ both to $T_{25}$ and to $V_{25}$.

Suppose for now that $y$ is connected to $V_{25}$ (we call this event
$E'_V$).
 Let $F'_V$ be the event that there is a vertex of $T_{25}$ 
 in the thin annulus $A_{|y-x|,|y-x|+0.05} $.
If $F'_V$ occurs,
then if $y$ comes from the first half-intensity process
$\Po_{\mu/2,24,25}$  
and there are no further vertices from 
$\Po_{\mu/2,24,25}$  (an event
of probability at least $\epsilon_4$),
event $F_V$ occurs.

If $F'_V$ does not occur, then let $V_{25}^y$ denote
the set of points of $\Po_{\mu,25,29}\cup \{y\}$ that are path-connected
to $V$, and build 
a process of intensity $\mu p(\cdot,T_{25})(1-p(\cdot, V_{25}^{y}))$,
  inwards inside $C_{|y-x|}$, until we get a vertex $z \in A_{24,|y-x|}$
(this must happen if $E'_V \cap 
H(V,T,\bS)$ is to occur but $F'_V$ does not occur). 
If then $y$ and $z$ both come from $\Po_{\mu/2,24,25}$
and there are no further vertices in $\Po_{\mu/2,24,25}$
(an event of probabilty at least $\epsilon_4$),
 then
$G_V$ occurs. Combining these yields
\begin{eqnarray}
P[F_V] + P[G_V]
 \geq
\epsilon_4
( P[E'_V \cap F'_V] +
 P[E'_V \cap H(V,T,\bS) \setminus F'_V])
\nonumber
\\
\geq
\epsilon_4
 P[E'_V \cap H(V,T,\bS)]. 
\nonumber
\end{eqnarray}
If $E' \setminus E'_V$ occurs, then $y$ is connected to $T_{25}$ and 
a simlar argument yields
$$
P[F_T] + P[G_T] \geq \epsilon_4 P[(E' \setminus E'_V) \cap H(V,T,\bS)], 
$$
and combining the last two estimates
 gives us (\ref{0507a}). \hfill $\Box$ \\

The following result, combined with Lemma \ref{Hcase1}, completes
the proof of Lemma \ref{Hprop}. 
\begin{lemm}
\label{Hcase2}
Suppose that $f$ is nonincreasing and (\ref{rangeone}) holds, and
 that $\rho \geq \frac{1}{\sqrt{2}} - 0.05$. Then the conclusion
of Lemma \ref{Hprop} holds. 
\end{lemm}

{\bf Proof.}
If $F_V$ or $F_T$  occurs we can continue in
similar fashion to the argument for Gilbert's  graph, as follows. 
Suppose $F_V$ occurs, let $y$ be as in the definition
of $F_V$ and set $r=|y-x|$,
 and let $z$ be an arbitrarily chosen point of
$T_{25}$ lying in $A_{r,r+.05}$.

\begin{figure}[htbp]
\includegraphics[angle = 0, width = 12cm]{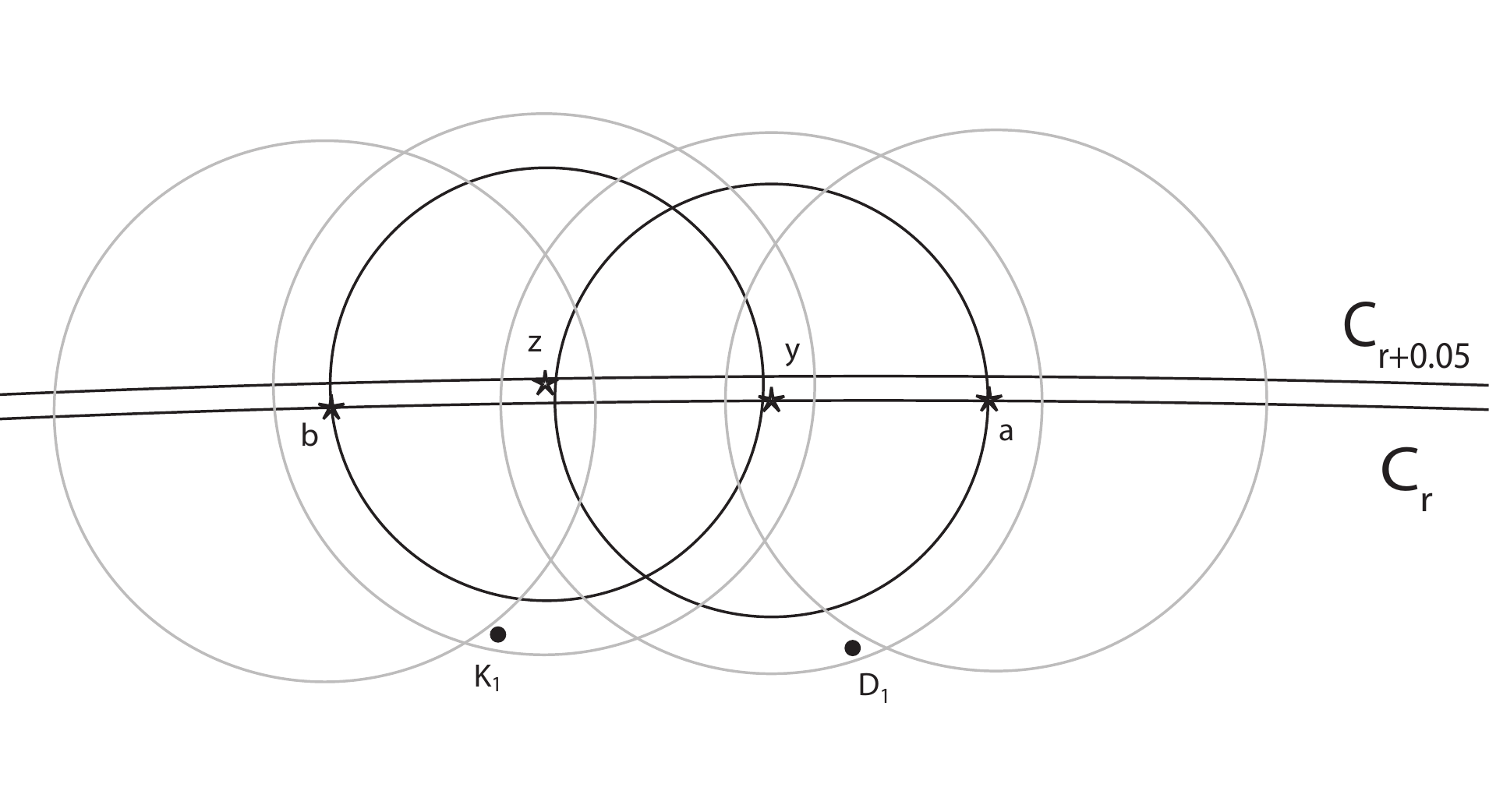}
 \caption{The grey circles are of radius $1$ and the black circles are
 of radius $\rho$.} \label{fig:points3}
\end{figure}

Let $a$ be the point with $|a-y| = \rho$  and
 $|a-x| = r$, on the other side of $y$
to $z$ (see Figure $\ref{fig:points3}$). 
Let $b$ be the point with
 $|b-z| = \rho$ and $|b-x|=r$, lying 
on the other side of $z$ to $y$. Let $a_1$ be the point
in $C_r$
with $|a_1-y| = 0.99$  and $|a_1-a| = 1.01$ 
 and let
$D_1:= C_{0.001}(a_1)$. Similarly let
$b_1 $ be the point in $C_r$ with $|b_1 -z| = 0.99$ and
$|b_1 - b| = 1.01$, and let $K_1:= C_{0.001}(b_1)$.
Note that $|y - z| > \rho$ so ${\rm dist}(D_1 , K_1) > \rho + 0.01$.

Let $D_2,\ldots, D_7$ and $K_2,\ldots, K_7$ be
 discs of radius $0.001$ and successive centres distant
$0.99$ from each other, such that, as before, having exactly
one red vertex in each of these little circles and no other vertices in
$A_{20,25}$, and connections between  the vertices in
successive circles $D_i,D_{i+1}$ and $K_j,K_{j+1}$
 ensures that $H'(V,T,\bS)$ occurs.

 Now sample $\Po'_{\mu/2,24,25} \cup \Po_{\mu,20,24}$ 
and  consider the event $E_5$, that there is exactly one new vertex $y_i$
in $D_i$ and exactly one new vertex $z_i$ in $K_i$ for $1 \leq i \leq 7$,
and no other new vertices. Then
$$
P[E_5 | F_V] \geq 
(0.001^2\pi\mu )^{14}
 \exp (-25^2 \pi \mu)=: \epsilon_5.
$$
Next, decide which edges are created from the new vertices.
 We want $y_1$ to connect
with $y$ (which happens with probability at least $f(0.991)$) but
not to any vertices in $T_{25}^{y,z}$ (which cannot happen as $a$ is the
closest place for a vertex in $T_{25}^{y,z}$). Similarly we also want $z_1$ to
connect with $z$ but not to $V_{25}^{y,z}$. 
Also we  want $z_1,y_1$ to not to be joined, and
we want connections betwen vertices in successive circles $D_i,D_{i+1}$
and $K_i,K_{i+1}$.  Given $F_V \cap E_5$,
these events all happen 
 with probability at least
$f(0.991)^{14} (1- f(\rho+0.01))$, in which case
$H'(V,T,S)$ occurs; hence
\begin{eqnarray}
P[H'(V,T,\bS)  |F_V]  \geq
f(0.991)^{14} [1-f(\rho+0.01)] \epsilon_5 := \epsilon_6
\label{0512b}
\end{eqnarray}
and similarly  $P[H'(V,T,\bS)|F_T] \geq \epsilon_6$.

Now suppose the event $G_V$ occurs. Then with $r:=|y-x|$,
 we have $z$ inside $C_{r}$
connected to a vertex $z_0$ of $T_{25}$ which must be outside $C_{r+0.05}$.

\begin{figure}[htbp]
\includegraphics[angle = 0, width = 12cm]{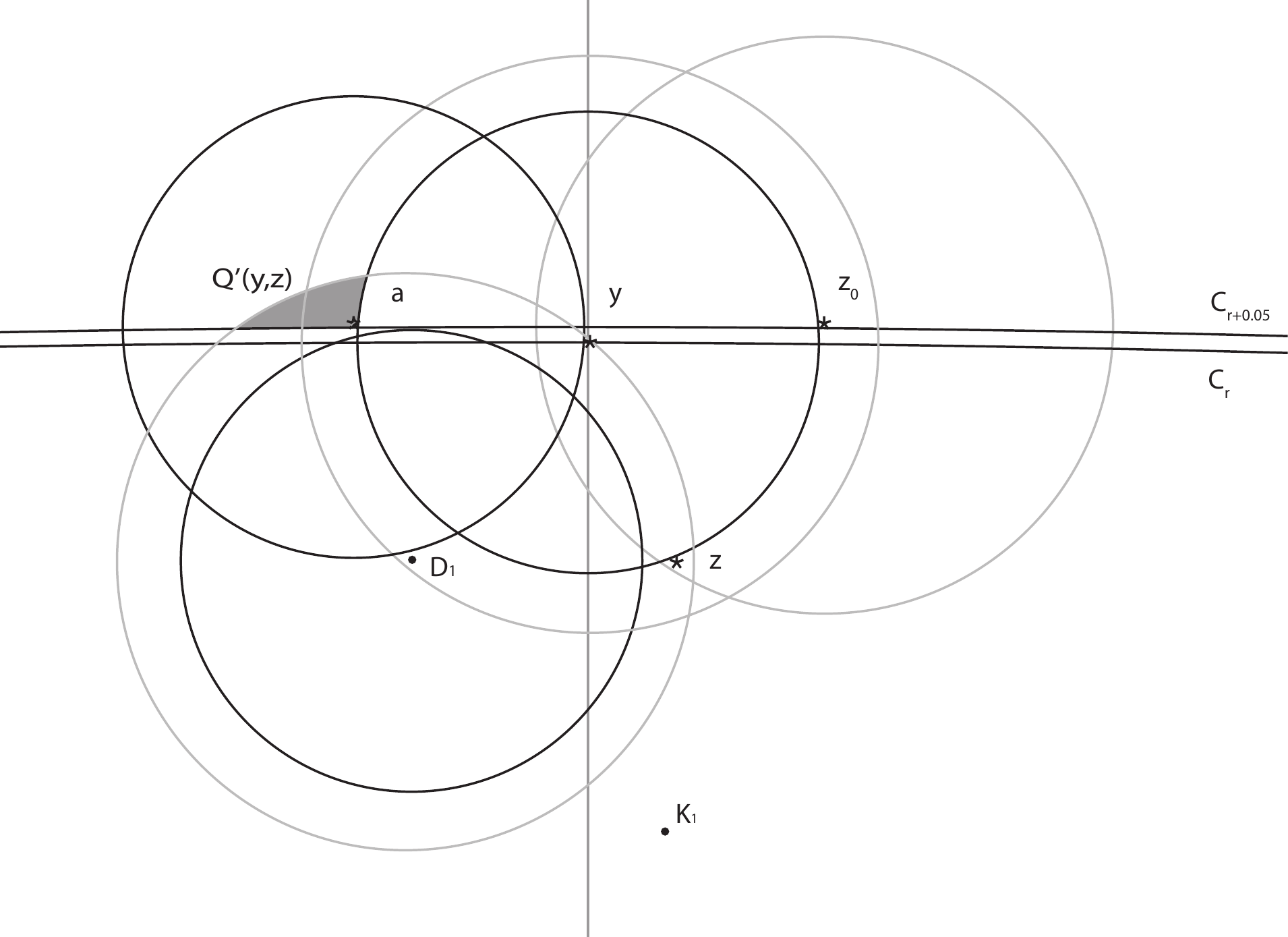}
 \caption{The grey circles are of radius $1.0001$ and the black circles are of
 radius $\rho$. }
 \label{fig:points4}
\end{figure}

 Let $a$ be the point with $|a=y|= \rho$  and
$|a-x|=r+0.05$ on the opposite side of $y$ to $z_0$. Let $l_a$ be the arc of $C_{r+0.05}$ to the left of $a$
 (see Figure $\ref{fig:points4}$). 
%
Let $b'(y,z)$ be
the point  at distance $0.999$ from $y$ and $\rho+0.01$ from $l_a$,
and define the region
$$
Q'(y,z):= C_{1.0001}(b'(y,z)) \setminus (C_{\rho}(y) \cup C_{r+0.05}(x)).
$$
The diameter of $Q'(y,z)$ is less than $0.9$, due
 to $\rho$ being at least $2^{-1/2} - 0.05$. Let
$E_6$ be the event that there are no more than $2$ vertices of 
$T_{25}^{y,z}$ in 
$Q'(y,z)$,
and no other vertices than those of $T^{y,z}_{24} \cup V^{y,z}_{24}$ in
 $A_{25,29}$.
By a similar argument to the proof of  Lemma $\ref{Q}$ the
 conditional probability of $E_6$ satisfies
\begin{eqnarray}
P[E_6 | G_V] \geq f(0.9)  \exp(-29^2 \mu) =: \epsilon_7.
\label{0512c}
\end{eqnarray}

Set $D_1 = C_{0.0001}(b'(y,z))$.  If there
is a vertex in $D_1$ it will be distant at least 
$\rho+0.001$ from $z$ and from any
vertex in $T^{y,z}_{25}$
 (as $l_a$ is the closest place such a vertex can be given
$G_V$ occurs) and at most $0.9991$ from $y$.
Let $a_1$ be the point
distant $0.999$ from $z$ on the line parallel with $xy$
through $z$, and let $K_1 := C_{0.0001}(a_1)$.
We can then pick little discs $D_i,K_i$, $2 \leq i \leq 7$, of
radius $0.0001$, as before
(see Figure \ref{fig:points1})
 such that if there is exactly one vertex in each of these discs and
 no other vertices
in the rest of $\Po_{20,25}$, and connections between vertices
in successive discs, then $H'(V,T,S)$ occurs.

Suppose that for the Poisson  process
$\Po_{\mu,20,24} \cup \Po'_{\mu/2,24,25}$
there is exactly one vertex
$y_i \in D_i$ and exactly one vertex $z_i \in K_i$ for each $i$
and there are no other vertices.
This has probability
at least
$[0.0001^2\pi \mu /2  ]^{14}
 \exp (-25^2 \pi \mu)$.
Given this event, consider now the event
 that
we get all connections occurring along the paths
$(y,y_1,\ldots,y_7)$ and $(z,z_1,\ldots,z_7)$
but no connection from $y_1$ to any vertex in $T_{25}^{y,z}$.  
 This has probability  at least
$ f(0.9991)^{14} [1-f(\rho+0.001)]^3$ (assuming $E_6$ occurs),
 and if this occurs
then $H'(V,T,\bS)$ occurs.
Hence,
$
P[ H'(V,T,\bS)| E_6 \cap G_V] \geq \epsilon_8$, 
 with 
$$
\epsilon_8 :=
[0.0001^2\pi \mu /2  ]^{14}
 \exp (-25^2 \pi \mu)
 f(0.9991)^{14} [1-f(\rho+0.001)]^3 \leq \epsilon_6.
$$

Combined with (\ref{0512b}) and (\ref{0512c}),
and a similar argument in the case of $G_T$,
this gives us 
(for $\rho \geq \frac{1}{\sqrt{2}} - 0.05$) the bound 
$$
4  P[H'(V,T,\bS) ] \geq \epsilon_7 \epsilon_8 (
P[G_V]+
P[G_T]+
P[F_V]+
P[F_T]).
$$
Combined with  Lemma \ref{lem0507}, this gives us
the desired result (\ref{Hpropeq}) with
$\epsilon = 0.25 \epsilon_4  \epsilon_7 \epsilon_8$.
%
\hfill{$\Box$} \\


\section{Proof of Theorem \ref{RCM}}
\eqco

We now generalise Theorem \ref{bool}  to
 the random connection model with connection
function $f:[0,\infty) \rightarrow [0,1]$, where $f$ is nonincreasing and has bounded support. 
Without loss of generality we assume (\ref{rangeone})
holds
 (as if not
we can rescale). For the enhancement this time we say that a vertex
$x$  is correctly configured if it is closed and 
 there are only $4$ vertices 
$v,w,y,z$ within $1$ of $x$, they are all red and joined to
$x$ and $v \sim w$ and $y \sim z$ but there are no other edges
amongst $v,w,y,z$. 
Notice that another
vertex could be not joined to $x$ but still cause it to be
incorrectly configured by being within $1$ of it.

All parts of the proof for this model are the same apart from Lemma 
\ref{intermediate}.
Accordingly we give a proof of the
 equivalent of Lemma \ref{intermediate} for the
Random Connection Model under our current assumptions.
\begin{lemm}
\label{RCMlemma}
Suppose $f$ is nonincreasing and (\ref{rangeone}) holds.
Then there is a continuous function $\delta :(0,1)^2 \to (0,\infty)$
such that for all $(p,q) \in (0,1)^2$, $n > 100$ and
$x \in B_n$,
\[
P_{n,2}(x) > \delta(p,q)P_{n,1}(x).
\]
\end{lemm}
In the proof we again write $C_r$ for $C_r(x)$.
Also we define events $E_{n,1}(x)$ and $R_n(x,\alpha,\beta)$
as in Section \ref{secgil}.
 It can easily be seen that the
proof of Lemma $\ref{intermed2}$ extends to this case as again the number of
possible green vertices in the completed process in a bounded
region is bounded.
Therefore
\begin{equation}\label{Heqn}
P[E_{n,1}(x) \cap R_n(x,20,30)] \geq \delta_1 P_{n,1}(x).
\end{equation}
Assume for now that $30.5 < |x| < n-30.5$.
Now suppose we create the whole process of intensity $\lambda$
in $B_n \setminus C_{30}$ and the red process of intensity $p \lambda$ in the
annulus $A_{29,30}$. We decide which vertices
outside $C_{30}$ are red, and assuming no closed vertices
occur in $A_{29,30}$, we then know which vertices 
outside $C_{30}$ are correctly configured.
We then determine which of these are green.

At this stage,
let $V$ be the set of coloured vertices in $B_n \setminus C_{29}$
that are connected (by a coloured path) to $B_{0.5}$ and
 let $T$ be the  coloured vertices in $B_n \setminus C_{29} $
that  are connected  to $\partial B_n$. Let $S$ be the remaining
coloured vertices in $B_n \setminus C_{29}$, and let
${\cal E}$ be the set of edges on $S$ inherited from the original random
connection model. Set $\bS := (S,{\cal E})$.

Then we can apply Lemma \ref{Hprop},
 using these boundary conditions, to the Poisson process
of red vertices,
of intensity $\mu = \lambda p$ inside $C_{29}$. 
If 
$E_{n,1}(x) \cap R_n(x,20,30)$ occurs, then $H(V,T,\bS)$
must occur, and therefore by Lemma \ref{Hprop},
$$
P[H'(V,T,\bS)]  \geq \epsilon(\lambda p) \delta_1 P_{n,1}(x).
$$
Now we can find $\delta_2$ such that
given $H'(V,T,\bS)$ occurs, 
 the probability of $x$ being $2$-pivotal is at least $\delta_2$. Indeed, with $y^*$ and $z^*$ as in the definition
of $H'(V,T,S)$, we just find little discs
$D_1,\ldots,D_{30}$ and $K_1,\ldots,K_{30}$
of radius $0.005$ leading from $y^*$ and $z^*$ in towards a bow-tie
 configuration around $x$ such that having one red vertex
 in each of these discs, with connections between successive discs,
 and no other vertices inside $C_{20}$, no
 vertices in the non-red process inside $C_{30}$, and having 
$Y_0 > p$ ensures $x$ is $2$-pivotal. This all occurs with probability
 at least
\[
\delta_2 := (0.005^2 \pi \lambda p )^{60} [f(0.9)]^{64} \exp (-900 \pi \lambda)(1-p). 
\]

Therefore we have 
\[
P_{n,2}(x) \geq  \delta_1 \delta_2 \epsilon P_{n,1}(x)
\]
for $30.5 < |x| < n-30.5$.

If $|x| \leq 30.5$ or $|x| \geq n- 30.5$,  then by minor modifications
of  the last part of the proof of Lemma \ref{intermediate}  we can
find  some  continuous $\delta_3:(0,1)^2 \to (0,\infty)$ such that
$P_{n,2}(x) \geq \delta_3 (p,q) P_{n,1}(x)$.
 So taking $\delta = \delta_1 \delta_2 \delta_3 \epsilon$ will give us the result. 
 \hfill{$\Box$} 

\section{Proof of Theorems \ref{th:squash} and \ref{th:spread}}
\eqco

For proving Theorem \ref{th:squash}, it is useful to
 consider mixed site-bond percolation on the graph $RCM(\lambda,f)$. Each site is open with probability $p$, and each bond is open with probability $q$. Clearly
the graph resulting from performing this mixed percolation process
on $RCM(\lambda,f)$ may be viewed as a  realisation of
$RCM(p\lambda, qf)$.

\begin{figure}[htbp]
\includegraphics[angle = 0, width = 12cm]{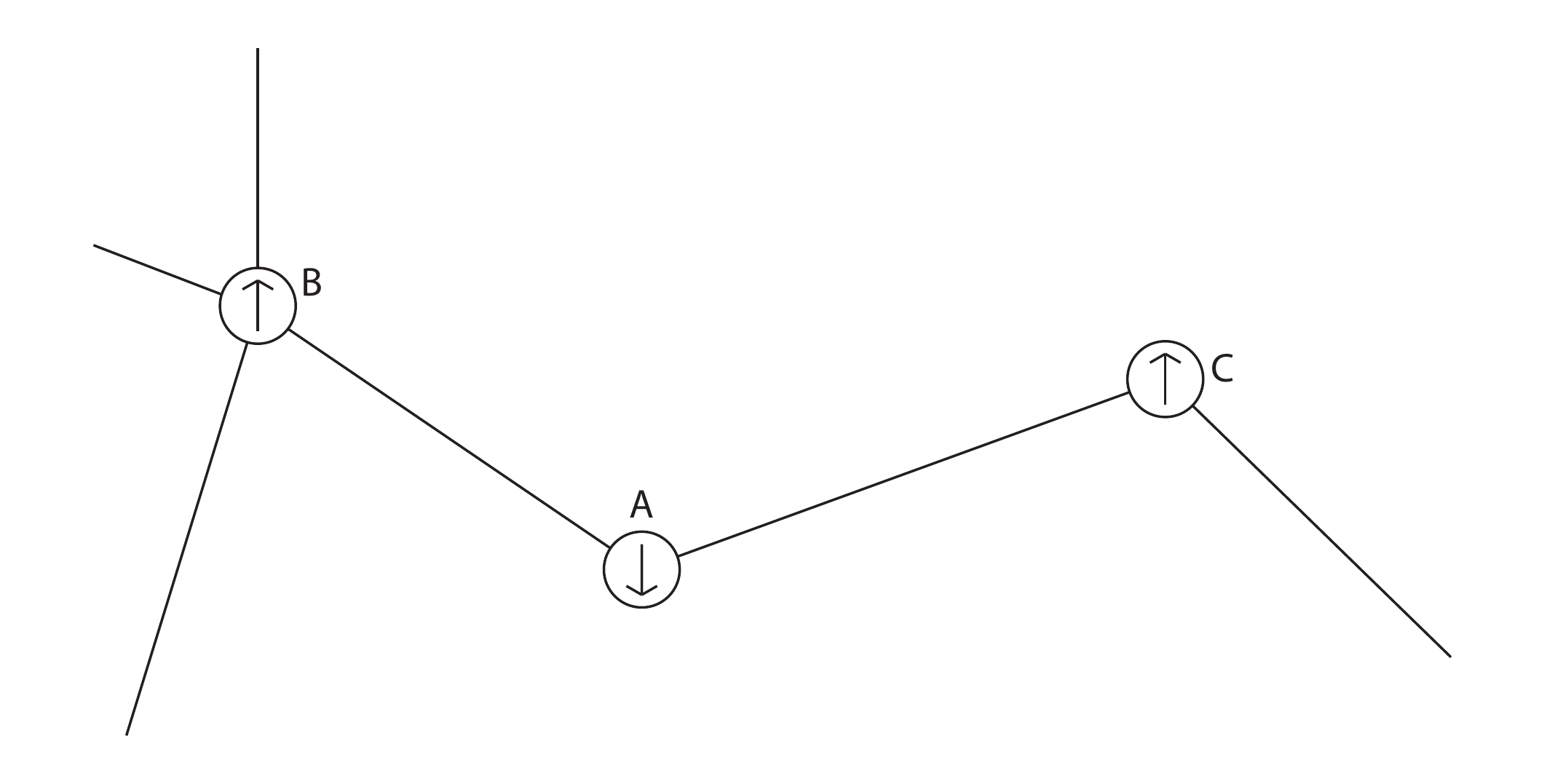}
 \caption{Here the vertex $A$ is correctly configured but $B$ and $C$ are not.} \label{fig:thm4}
\end{figure}

In proving Theorem \ref{th:squash} we assume without loss of 
generality that (\ref{rangeone}) holds.
  We consider a new site percolation model, where sites are open with probability $pq$ if they are correctly configured and with probability $p$ if they are not correctly configured. 
Each site is designated either an {\em up-site} or a {\em down-site},
each with probability $1/2$, independently of everything else.
We say vertex $y$ is a 1-neighbour of vertex $x$ if $|y-x| \leq 1$.
A vertex $x$ is correctly configured if 
it has exactly two 1-neighbours
 (denoted $y_1$ and $y_2$, say)
and $x$ is connected both to $y_1$
and to $y_2$,  and $x$ is a down-site but $y_1$ and $y_2$ are up-sites
(see Figure \ref{fig:thm4}).

The extra randomization of up-sites and down-sites is designed
 to ensure that if a site is correctly
configured, then its neighbours are not.

 We build this model by having a Poisson process of intensity $\lambda$ and labelling vertices $x_1, x_2,\ldots$ in order of distance from the origin.
 We also have  independent uniform random variables $W_i,Y_i, Z_i$ for
$i=0,1,2,\ldots$.
We say vertex $x_i$ is  up-site if and only if $W_i < 1/2$.
 If a vertex $x_i$ is correctly configured it is open if
 $Y_i < p$ and $Z_i < q$. Otherwise it is open if $Y_i < p$. We define $\partial B_n$ to be $B_n \setminus B_{n-0.2}$. We let $A_n$ be the event that there is an open path from $B_{0.2}$ to $\partial B_n$ in the process restricted
 to $B_n$, and for $x \in B_n$ define
 $A_n^x$ similarly in terms of the process in $B_n$
with an added vertex at $x$. 

Let the sites $x_i$ for which $Y_i <p$ be denoted {\em red} (a
standard Bernoulli site percolation process).
The set of open sites
may be viewed as a {\em diminishment} of the set of red sites,
in which each correctly configured red site is removed 
with probability $1-q$. We can couple the diminished site
percolation process to the mixed site-bond process (with parameters
$p,q$) in such a way that if the the mixed process percolates
then so does the diminished site process,
as follows.

List the edges of this graph in an arbitrary order as
$e_1,e_2,\ldots$,  and determine the open sites
and edges for the mixed site-bond process.
Deem  each vertex  to be red  if and only if it is 
open in the mixed process.  If a vertex $x$ is correctly configured, 
then it has degree  2 and 
  and has no correctly configured neighbour.
In this case, let $x$ be diminished (i.e. removed
from the set of red vertices) 
if and only if the first edge incident to it 
(according to the given ordering) is closed.

If there is an infinite open path in the mixed
percolation process,  we can find such a path
which starts at a non-correctly configured 
vertex. In this case, every vertex in
the path will be red and undiminished, so there
will be an infinite path in the diminished site 
percolation process as well.

Let $\theta_n(p,q)$ be the probability that $A_n$ occurs and 
 let $\theta(p,q)$ be the limit inferior.
 The proof of Proposition $\ref {prop:enh}$ is easily modified to this model.  We say that a point $x$ is $3$-pivotal if putting a vertex at $x$ and making $Y_0 <p$ means that $A^x_n$ occurs but having $Y_0 > p$ means its does not. Similarly with $4$-pivotal and $Z_0$. Again we have a form of the Margulis-Russo formulae:
\bea
\frac{\partial \theta_n(p,q)}{\partial p}=\int^{}_{B_n} \lambda P_{n,3}(x,p,q) \, \mathrm{d}x
\label{eq11}
\eea
and
\bea
\frac{\partial \theta_n(p,q)}{\partial q} = \int^{}_{B_n} \lambda P_{n,4}(x,p,q) \, \mathrm{d}x.
\label{eq12}
\eea

We then need to prove the equivalent of Lemma \ref{intermediate}:
\begin{lemm}
\label{lem1004}
There is a continuous function $\delta:(0,1)^2 \to (0,\infty)$ such that
for all $n > 100$ and all $x \in B_n$ we have
\begin{equation}
 P_{n,4}(x,p,q) \geq \delta (p,q) P_{n,3}(x,p,q).
\label{0420a}
\end{equation}
\end{lemm}

Before proving this, 
we give the equivalent of Lemma \ref{intermed2} which
says we can assume all the vertices in an annulus of fixed
size are red and none of them are diminished.
Given $p$ and $q$, and given $2 < \alpha < \beta$,
 let $R_{n}(x,\alpha,\beta)$  be the event that 
all vertices in $A_{\alpha,\beta}(x) $ are red.
Let 
  $R'_{n}(x,\alpha,\beta)$  be the event that 
  $R_{n}(x,\alpha,\beta)$  occurs and also none of
the vertices in $A_{\alpha-1,\beta+1}(x)$  is diminished.
 We claim there exists continuous $\delta_1:(0,1)^2 \to (0,\infty)$ such
that for all  $n$ and $x$ we have
\begin{equation}
\label{0420b}
 P[E_{n,3}(x) \cap R'_n(x,\alpha,\beta)]  \geq \delta_1(p,q) P_{n,3}(x,p,q).
\end{equation}
To prove this we create the whole Poisson point  process of intensity $\lambda$
 in $B_n$, and the edges between these vertices, and
decide which vertices 
 outside $A_{\alpha-1,\beta +1}$ are red, and which of them are up-sites.
For each vertex in $A_{\alpha-1,\beta +1}$
  having more than two 1-neighbours and/or having a down-site outsite
 $A_{\alpha-1,\beta +1}$ as a 1-neighbour, we decide if that vertex
is red, and whether it is an up-site or a down-site (these vertices
cannot be correctly configured). We then know which of the vertices outside
$A_{\alpha-1,\beta-2}$ are correctly configured and we decide
which of them are open.

 This leaves a set $W $ of vertices
in $A_{\alpha-1,\beta +1}$ with at most two neighbours which are the ones that could be correctly
 configured. As at (\ref{0513a}), the set $W$ has at most $12 (\beta +2)^2$
 elements and for $x$ to have a chance of being 
 3-pivotal there must exist a subset
 $W'$ of $W$ such
 that if all the vertices in $W'$ are open and all the vertices in
 $W \setminus W'$ are closed then $x$ is 3-pivotal.
So if $Y_i < p$ for all $x_i$ in
 $W'$ and $Y_i > p$ for all $x_i$ in $W \setminus W'$,
 then the event 
$E'_{n,3}(x)$ occurs, where
$E'_{n,3}(x)$ denotes the event that $x$ is 3-pivotal
in a modified model where the diminishments are suppressed
in $A_{\alpha-1,\beta +1}(x)$. Hence
$$
P[E'_{n,3}(x) ] \geq [p(1-p)]^{12(\beta+2)^2}P_{n,3}(x).  
$$
Adding or removing extra non-red vertices in $A_{\alpha,\beta}(x)$
does not affect event $E'_{n,3}(x)$ and therefore $E'_{n,3}(x)$
is  independent of the event $R_{n}(x,\alpha,\beta)$.
Also, if $E'_{n,3}(x) \cap R_n(x,\alpha,\beta)$ occurs, then there
are at most $12(\beta +2)^2$ correctly configured
red vertices 
in $A_{\alpha-1,\beta +1}$,
 and the probability
that none of these is diminished is at least $(1-q)^{12(\beta+2)^2}$.
In this case event $E_{n,3} \cap R'_n(x)$ occurs, and (\ref{0420b}) follows with
$$
\delta_1 := [p(1-p)(1-q)]^{12 (\beta+2)^2}  
\exp(- \pi (\beta^2 - \alpha^2) \lambda (1-p)).
$$  

{\bf Proof of Lemma \ref{lem1004}.}  
 Assume for now that $30.5 < |x| < n- 30.5$.
Create the full process of intensity $\lambda$ outside
 the circle $C_{30}$ around $x$, and
the red process $\Po_{p \lambda,29,30}$ in the annulus $A_{29,30}(x)$.
 Assuming there are no other vertices
in $A_{29,30}(x)$, determine which
 vertices outside $C_{30}$ are diminished, but do not
yet diminish any vertices inside $C_{30}$. Deem open all
vertices that are red and have not been diminished 
at this stage. 
 Let $V$ be the set of vertices
 now connected (by an open path) to $B_{0.2}$ and let $T$ 
be those vertices
 connected to $\partial B_n$.
Let $S$ be the remaining open  vertices,
and let ${\cal E}$ be the edges on $S$ inherited
the original random connection model. Set $\bS =(S,{\cal E})$. 

Now create the red process  $\Po_{p \lambda,20, 29}$.
 Let events $H := H(V,T,\bS)$ and $H' := H'(V,T,\bS)$ 
 be as described just before Lemma \ref{Hprop} (with $\mu = p \lambda$).

Event $H$ must happen if $E_{n,3}(x) \cap R'_n(x,20,29)$
is to occur.
Hence by (\ref{0420b}),
$P[H] \geq \delta_1 P_{n,3}(x)$, and therefore by
Lemma \ref{Hprop}, 
 $P[H'] \geq  \epsilon(p \lambda ) \delta_1 P_{n,3}(x)$. 

As in the latter part of the proof of Lemma \ref{RCMlemma},
if $H'$ occurs 
we can then form little discs $D_1,\ldots,D_{30}$ forming
a path in $C_{20}(x)$ from $y^*$ to $x$ and 
$K_1,\ldots,K_{30}$ forming a path in $C_{20}(x)$ from $z^*$  to
 $x$, this time with only $D_{30}$ and $K_{30}$ within
unit distance of $x$. Then 
$x$ will be $4$-pivotal
if we have exactly
 one red vertex in each of these discs,  no other vertices in the rest of 
the process $ \Po_\lambda \cap C_{30}$, 
  all edges along these paths are present,
 $Y_0 < p$, $x$ is a down-site but its neighbours are up-sites, and
no vertices in $C_{30}$ are diminished. 
 This all occurs with probability at least 
\[
\delta_2(p,q) 
:= (0.005^2 \pi \lambda p )^{60} [f(0.9)]^{62} \exp (-900 \pi \lambda)(p/8)
(1-q)^{12(31^2)}, 
\]
where the last factor is a lower bound on the probability
that no diminishment occurs in $C_{30}$, by
the same argument as in the proof of (\ref{0420b}).
Therefore
\[
P_{n,4}(x) \geq \delta_1 \delta_2  \epsilon(p \lambda) P_{n,3}(x),
\]
for $x$ with $30.5 < |x| < n-30.5$.
 For other $x$ we can argue similarly
 to the last part of the proof of Lemma \ref{intermediate}  to
find  some  continuous $\delta_3:(0,1)^2 \to (0,\infty)$ such that
$P_{n,4}(x) \geq \delta_3 (p,q) P_{n,3}(x)$.
 So taking
$
\delta = \min(  \delta_1\delta_ 2 \epsilon,\delta_3),
$
 we are done.
\hfill{$\Box$} \\

{\bf Proof of Theorem \ref{th:squash}.}
We take $q_0 < 1$, and fix $f$.
We now set $q^* = (1+q_0)/2$, and
choose $\lambda > \lambda_{q_0f}$,
and consider the graph  $RCM(\lambda,f)$.
 Define $p^{*} := \lambda_{q_0 f}/\lambda$, so 
$p^* \in (0,1) $.
Now by considering a small box around $(p^*,q^*)$, and using
Lemma \ref{lem1004}, \eq{eq11}, \eq{eq12} and the analogue of
Proposition \ref{prop:enh},
we can find  $\epsilon > 0 $
such that $\epsilon < \min(p^*,1-p^*)$ and
\[
\theta(p^* + \epsilon,q_0) \leq \theta (p^* ,q^*) 
\leq \theta(p^* - \epsilon,1).
\]
Now the definition of $p^*$ implies that
 $RCM( (p^* + \epsilon)\lambda,q_0f)$ percolates.
Hence, on $RCM(\lambda,f)$ the mixed site-bond process with parameters
$(p^*+\epsilon,q_0)$ percolates and therefore
the diminished site process with parameters $(p^*+\epsilon,q_0)$
percolates.  Thus
  $\theta (p^* + \epsilon ,q_0) > 0$ and hence
  $\theta(p^* - \epsilon,1 ) > 0$
 which means that $RCM((p^* - \epsilon)\lambda,f)$ percolates so
\[
\lambda_f \leq (p^* - \epsilon)\lambda < p^* \lambda = \lambda_{q_0 f}
\]
and we are done.  \hfill{$\Box$} \\

{\bf Proof of Theorem \ref{th:spread}.}
Since $S_pf \equiv S_{p/q} S_qf$, it suffices to consider the
case with $q=1$ so that $S_qf \equiv f$.
Define the connection function
$g(r) = f(\sqrt{p} r)$. Then $S_pf \equiv p g$ so by
Theorem \ref{th:squash},
$$
\lambda_{S_pf} > \lambda_{g}
$$
and hence the graph RCM$(\lambda_{S_pf},g )$ 
is the realization of a supercritical random connection
model.  Let $p_c^{\rm bond} $ and $p_c^{\rm site}$
denote the critical values for bond, respectively site,
percolation on this graph.
By Theorem \ref{RCM}, we have
$p_c^{\rm bond} > p_c^{\rm site}$.

Given $p' \in(0,p)$, we have $p'g \equiv (p'/p) S_pf $
 so that $\lambda_{S_pf}  < \lambda_{p'g}$ by Theorem
\ref{th:squash}, so that the graph $RCM(\lambda_{S_pf},p'g)$
does  not percolate, and therefore $p_c^{\rm bond} \geq p'$.
Hence, taking $p' \uparrow p$, we have
\begin{eqnarray}
p \leq p_c^{\rm bond} < p_c^{\rm site}.
\label{MP0224}
\end{eqnarray}
By scaling, the graph RCM$(\lambda_{S_pf},g)$ is
equivalent to the graph RCM$(p^{-1} \lambda_{S_pf},f )$
so that 
$$
p_c^{\rm site} = \frac{p \lambda_f}{\lambda_{S_pf}}
$$
and combining  this with (\ref{MP0224})
yields
the desired inequality
 $\lambda_{S_pf} < \lambda_{f}$.

For the last part, observe that
 $\int_0^\infty r S_pf(r) dr = \int_0^\infty r f(r) dr$ and therefore 
the fact that 
(\ref{100122a})
holds as  a weak inequality 
 for $S_pf$ implies that it  holds as a strict inequality
for $f$.  \hfill{$\Box$}


\end{document}